\documentstyle[12pt]{amsart}

\renewcommand{\qed}{\hfill{\em q.e.d.}}
\newcommand\codim{{\operatorname{codim}}}

\newcommand\Pic{{\operatorname{Pic}}}
\newcommand\Div{{\operatorname{Div}}}

\newcommand\supp{{\operatorname{supp}}}
\newcommand\pdiv{{\operatorname{div}}}

\newcommand\depth{{\operatorname{depth}}}
\newcommand{\prs}[1]{{\Bbb P}^{#1}}
\newcommand{\hyq}[1]{{\mbox{\boldmath $Q$}^{#1}}}

\newcommand{\lsys}[1]{{\mid\;#1\;\mid}}
\newcommand{\str}{{{\cal O}}}

\newcommand{\nbl}[1]{{{\cal N}_{#1}}}

\newcommand{\dsh}[1]{{\omega^0_{#1}}}
\newcommand\euch[2]{{\chi(#1,#2)}}

\newcommand\rational{{\Bbb Q}}
\newcommand\real{{\Bbb R}}
\newcommand\integer{{\Bbb Z}}
\newcommand\complex{{\Bbb C}}

\newcommand\sing{{\mbox{sing}}}

\newcommand\non{{\mbox{non}}}

\newcommand\Etil{{\tilde{E}}}

\newcommand\minus{{\setminus}}
\newcommand{\rest}[2]{{{#1}\!\mid_{#2}}}
\newcommand{\shortseq}[5]{{0\to{#1}\stackrel{#2}{\to}{#3}
\stackrel{#4}{\to}{#5}\to0}}
\newcommand{\Coh}[3]{{H^{#1} ({#2},{#3})}}
\newcommand{\coh}[3]{{h^{#1} ({#2},{#3})}}

\newtheorem{th}{Theorem}[section]
\newtheorem{de}[th]{Definition}
\newtheorem{ex}[th]{Example}
\newtheorem{lem}[th]{Lemma}
\newtheorem{lede}[th]{Definition and Lemma}
\newtheorem{co}[th]{Corollary}
\newtheorem{cl}[th]{Claim}
\newtheorem{re}[th]{Remark}
\newtheorem{pr}[th]{Proposition}
\newtheorem{con}[th]{Conjecture}
\newtheorem{fa}[th]{Fact}
\newtheorem{oq}[th]{Question}
\newtheorem{prob}[th]{Problem}
\newtheorem{as}[th]{Assumption}
\newenvironment{TH}{\begin{th}\em}{\end{th}}
\newenvironment{DE}{\begin{de}\em}{\end{de}}
\newenvironment{EX}{\begin{ex}\em}{\end{ex}}
\newenvironment{LE}{\begin{lem}\em}{\end{lem}}
\newenvironment{LEDE}{\begin{lede}\em}{\end{lede}}

\newenvironment{CL}{\begin{cl}\em}{\end{cl}}
\newenvironment{RE}{\begin{re}\em}{\end{re}}
\newenvironment{PR}{\begin{pr}\em}{\end{pr}}
\newenvironment{CON}{\begin{con}\em}{\end{con}}
\newenvironment{FA}{\begin{fa}\em}{\end{fa}}
\newenvironment{OQ}{\begin{oq}\em}{\end{oq}}
\newenvironment{PROB}{\begin{prob}\em}{\end{prob}}
\newenvironment{AS}{\begin{as}\em}{\end{as}}

\oddsidemargin=\evensidemargin
\title[On the fourth adjoint contractions]{On the fourth adjoint contractions\\
of divisorial and fiber types}

\author{Shu Gilbert Nakamura}

\curraddr{Department of Mathematics, University of California, Riverside, 92521-0001}

\email{ngilbert@@math.ucr.edu, snakamur@@kenna.math.nd.edu}

\begin{document}
\begin{abstract}
In this paper, we will list up all the cases for the ray contractions of divisorial and fiber types for smooth projective varieties of dimension {\em five}.
These are obtained as a corollary from the lists of 
$n$-dimensional $k$-th adjoint contractions
$f:X\to Y$
of the same types for
$k=1,2,3$
and
$4$
($n\geq5$).
The lists for
$k=1,2$
and
$3$
have previously been obtained in \cite[Proposition~1.2 and Theorem~1.3]{na:3adc}.
The main task will be to have such a list for
$k=4$, 
where one case in the list fails to show that a positive-dimensional general fiber
$F$
of
$f$
is irreducible when
$n>5$.
This assertion will, however, be proven when
$n=5$
with an essential aid of $3$-dimensional Minimal Model Program in \cite{mo:3mmp}.
(We do not show the existence of cases.)
\end{abstract}
\maketitle

\section{Introduction and main results}
\label{intro}
\noindent
Let 
$X$
be a 
$5$-dimensional smooth projective variety over the complex number field
$\complex$.
Assume that a canonical divisor
$K_X$
of
$X$
is not nef.
Then, Mori-Kawamata theory (\cite{kmm}, \cite{mo:3fd}) provides an {\em extremal ray}
$R$ of
$X$
which is defined by
$R=H^\bot \cap \overline{NE}(X),$
where
$H$
is a nef divisor of 
$X$
(called a {\em supporting divisor} of
$R$)
and 
$\overline{NE}(X)$
is the cone of curves of
$X$,
such that for any 
$\gamma\in R\minus\{0\}$,
$K_X \cdot \gamma<0$
and that
$\dim_{\real} R\otimes\real=1$.
Moreover, the linear system
$\lsys{mH}$,
$m>>0$,
gives rise to a morphism
$f:X\to Y$
onto a normal projective variety
$Y$
with only connected fibers such that
$f(C)$
is a point if and only if the numerical class
$[C]$
is in
$R$
(or equivalently,
$H\cdot C=0$)
for any integral curve in 
$X$,
and that 
$H=f^* {\cal A}$
for some ample Cartier divisor
${\cal A}$
of
$Y$.
We call (the isomorphism class of)
$f$ 
the {\em ray contraction} of 
$R$.
If the dimension of the exceptional locus of 
$f$
is equal to 
$\dim X-1$
or
$\dim X$,
$f$
is said to be {\em of divisorial type},
{\em of fiber type}, respectively.
The main result in this paper is the following
\begin{TH}
\label{frfd}
The tables {\bf A} and {\bf B} list up all the possible cases for the ray contractions
$f:X\to Y$ 
of divisorial and fiber types for smooth five folds
$X$. \\
{\bf A}.
Let
$f:X\to Y$
be of divisorial type.
Let
$F$
denote a general fiber of
$E\to f(E)$.
If
$\dim f(E)=1$
and
$F\cong\prs{3}$,
then
$F$ 
can be any fiber of
$E\to f(E)$
over a smooth point of
$f(E)$.
$$\begin{array}{||c|l|c||}\hline
\dim f(E)&\mbox{structure of a positive-dimensional fiber of $f$}&\mbox{singularities of }Y\\
\hline\hline
0&E\mbox{ is a Mukai variety and }\rest{K_X}{E}\cong\nbl{E/X}&\mbox{$1$-factorial, terminal}\\ 
\cline{2-3}
&E\mbox{ is a Del Pezzo variety and }2\rest{K_X}{E}\cong\nbl{E/X}&\mbox{$2$-factorial, terminal}\\
\cline{2-3}
&E\mbox{ is a Del Pezzo variety and }\rest{K_X}{E}\cong2\nbl{E/X}&\mbox{$1$-factorial, terminal}\\ 
\cline{2-3}
&E\cong\hyq{4},\rest{K_X}{E}\cong\str_{\hyq{}} (-1), \mbox{ and }\nbl{E/X}\cong\str_{\hyq{}} (-3)&\mbox{$3$-factorial, terminal}\\
\cline{2-3}
&E\cong\hyq{4},\rest{K_X}{E}\cong\nbl{E/X}\cong\str_{\hyq{}} (-2)&\mbox{$2$-factorial, terminal}\\
\cline{2-3}
&E\cong\hyq{4},\rest{K_X}{E}\cong\str_{\hyq{}} (-3), \mbox{ and }\nbl{E/X}\cong\str_{\hyq{}} (-1)&\mbox{$1$-factorial, terminal}\\
\cline{2-3}
&E\cong\prs{4},\rest{K_X}{E}\cong\str_{\prs{}} (-1), \mbox{ and }\nbl{E/X}\cong\str_{\prs{}} (-4)&\mbox{$4$-factorial, terminal}\\ 
\cline{2-3}
&E\cong\prs{4},\rest{K_X}{E}\cong\str_{\prs{}} (-2), \mbox{ and }\nbl{E/X}\cong\str_{\prs{}} (-3)&\mbox{$3$-factorial, terminal}\\ 
\cline{2-3}
&E\cong\prs{4},\rest{K_X}{E}\cong\str_{\prs{}} (-3), \mbox{ and }\nbl{E/X}\cong\str_{\prs{}} (-2)&\mbox{$2$-factorial, terminal}\\ 
\cline{2-3}
&E\cong\prs{4},\rest{K_X}{E}\cong\str_{\prs{}} (-4), \mbox{ and }\nbl{E/X}\cong\str_{\prs{}} (-1)&\mbox{smooth}\\ 
\hline
\end{array}$$

$$\begin{array}{||c|l|c||}\hline
1&F^3 \mbox{ is a Del Pezzo variety and }\rest{K_X}{F}\cong\rest{\nbl{E/X}}{F}&\mbox{$1$-factorial, terminal}\\ 
\cline{2-3}
&F\cong\hyq{3}, \rest{K_X}{F}\cong\str_{\hyq{}} (-1),\mbox{ and }\rest{\nbl{E/X}}{F}\cong\str_{\hyq{}} (-2)&\mbox{$2$-factorial, terminal}\\
\cline{2-3}
&F\cong\hyq{3}, \rest{K_X}{F}\cong\str_{\hyq{}} (-2),\mbox{ and }\rest{\nbl{E/X}}{F}\cong\str_{\hyq{}} (-1)&\mbox{$1$-factorial, terminal}\\
\cline{2-3}
&F\cong\prs{3}, \rest{K_X}{F}\cong\str_{\prs{}} (-1),\mbox{ and }\rest{\nbl{E/X}}{F}\cong\str_{\prs{}} (-3)&\mbox{$3$-factorial, terminal}\\
\cline{2-3}
&F\cong\prs{3}, \rest{K_X}{F}\cong\rest{\nbl{E/X}}{F}\cong\str_{\prs{}} (-2)&\mbox{$2$-factorial, terminal}\\
\cline{2-3}
&F\cong\prs{3}, \rest{K_X}{F}\cong\str_{\prs{}} (-3),\mbox{ and }\rest{\nbl{E/X}}{F}\cong\str_{\prs{}} (-1)&\mbox{$1$-factorial, terminal}\\
\hline
2&F\cong\hyq{2}, \rest{K_X}{F}\cong\str_{\hyq{}} (-1),\mbox{ and }\rest{\nbl{E/X}}{F}\cong\str_{\hyq{}} (-1)&\mbox{$1$-factorial, terminal}\\
\cline{2-3}
&F\cong\prs{2}, \rest{K_X}{F}\cong\str_{\prs{}} (-1),\mbox{ and }\rest{\nbl{E/X}}{F}\cong\str_{\prs{}} (-2)&\mbox{$2$-factorial, terminal}\\
\cline{2-3}
&F\cong\prs{2}, \rest{K_X}{F}\cong\str_{\prs{}} (-2),\mbox{ and }\rest{\nbl{E/X}}{F}\cong\str_{\prs{}} (-1)&\mbox{$1$-factorial, terminal}\\
\hline
3&F\cong\prs{1}, \rest{K_X}{F}\cong\rest{\nbl{E/X}}{F}\cong\str_{\prs{}} (-1)&\mbox{$1$-factorial, terminal}\\
\cline{2-3}
\hline
\end{array}$$

\smallskip

\noindent
{\bf B}.
Let
$f:X\to Y$
be of fiber type.
Let
$F$
denote a general fiber of
$f$.

$$\begin{array}{||c|l|c||}
\hline
\dim Y&\mbox{structure of a positive-dimensional fiber of $f$}&\mbox{singularities of }Y\\
\hline\hline
0&X\mbox{ is a Fano manifold of co-index five}&\mbox{smooth}\\
\cline{1-3}
1&F^4 \mbox{ is a Fano manifold of co-index four}&\mbox{smooth}\\ \cline{1-3}
2&F^3 \mbox{ is a Mukai manifold}&\mbox{$1$-factorial, normal}\\ \cline{1-3}
3&F^2 \mbox{ is a Del Pezzo manifold}&\mbox{$1$-factorial, normal}\\
\cline{1-3}
4&F\cong\hyq{1} \mbox{(smooth),  and }\rest{K_X}{F}\cong\str_{\hyq{}} (-1)&\mbox{$1$-factorial, normal}\\
\cline{1-3}
\hline
\end{array}$$
\noindent{where} 
$\hyq{}$ 
denotes a hyperquadric in a projective space.
For the definitions of a Fano variety, co-index, a Del Pezzo variety and a Mukai variety, see Definition~\ref{spva}.
(We do not show the existence of any case in the above table.)
\end{TH}

\medskip

\noindent
The way we see the problem is somewhat more general than that of the $5$-dimensional ray contractions.
This idea may be explained as follows:
Note that
$H-K_X$
is 
$f$-ample. 
Call this 
$L$,
i.e.,

$$\begin{array}{rcl}
H&=&K_X +1\cdot L\\
&=&K_X +(n-k)L,
\end{array}$$

\noindent
where
$n=\dim X=5$
and
$k=4$.
We can take
$L$
to be ample, and thus, our object becomes a polarized manifold
$(X,L)$
of dimension
$n$
with an adjoint divisor
$K_X +(n-k)L$,
where
$k=4$.
Furthermore, we do not like to have 
$L$
be ``numerically'' some multiple of another ample divisor
$A$,
since, for instance, if 
$L\approx2A$,
then

$$\begin{array}{rcl}
K_X +(n-k)L&\approx&K_X +2(n-k)A\\
&\approx&K_X +(n-k')A,
\end{array}$$

\noindent
where 
$k'<k$.
Hence, this case would essentially be reduced to an easier adjoint system:
$K_X +(n-k')A$.
This argument leads to a good interaction of both the Mori-Kawamata and the adjunction theories and produces the following proposition whose proof is given in Appendix.

\begin{PR}
\label{1stcl}
Let 
$X$ 
be an 
$n$-dimensional smooth projective variety over the complex number field 
$\Bbb C$ 
with the canonical divisor 
$K_X$ 
not nef
$(n\geq2)$.
Let
$R$
be an arbitrary extremal ray on 
$X$,
and let
$f:X\to Y$ 
be the ray contraction of 
$R$.
Then there are an ample Cartier divisor
$A$
on
$X$
and a unique integer
$-1\leq k\leq n-1$
such that
\begin{enumerate}
\item $R=(K_X +(n-k)A)^\bot \cap\overline{NE}(X)$
and that
\item $A$
is {\em numerically reduced} in the sense that there are no ample Cartier divisors
$A'$ 
satisfying
$A\cdot C=pA'\cdot C$
for some integer
$p>1$
and a nonzero effective curve
$C$
in
$R$.
\end{enumerate}
\end{PR}
\begin{DE}
\label{defadj}
In the setting of Proposition~\ref{1stcl}, we will call
$f:X\to Y$ 
the $k$-th {\em adjoint contraction} with a supporting divisor 
$K_X +(n-k)A$.
In particular,
$f$
is called the {\em fourth adjoint contraction} if
$k=4$.
\end{DE} 

\begin{RE}
We note that this definition has already divided all the ray contractions into
$(n+1)$-distinct families as
$k$
varies from
$-1$
to 
$n-1$.
\end{RE}

\noindent
All the
$k$-th adjoint contractions for 
$k$
ranging from
$-1$
to
$3$
are listed up in \cite{na:3adc}.
One wishes to extend the list to the case
$k=4$,
the fourth adjoint contractions.
As part of this attempt, we will obtain a list of all the fourth adjoint contractions of fiber and divisorial types, the main objects in this paper, in Theorem~\ref{fmth} below.
(Those of flipping type are not considered.)

\begin{TH}
\label{fmth}
Let
$f:X\to Y$
be an $n$-dimensional fourth adjoint contraction of fiber type or of divisorial type. 
(From Proposition~\ref{1stcl}, it follows that
$n\geq5$.)
Then, all the possible cases of
$f$
are listed up in the following table.
If
$f$
is of fiber type,
$F$
denotes a general fiber of
$f$.
If
$f$
is of divisorial type,
$E$
denotes a unique exceptional prime divisor of
$f$,
and
$F$
denotes a general fiber of
$\rest{f}{E}:E\to f(E)$
(or any fiber of
$\rest{f}{E}$
over a smooth point of
$f(E)$
if
$r=1$
and
$F\cong\prs{n-2}$),
where
$r=\dim f(E)$.

$$\begin{array}{||c|l|c||}\hline
\dim Y&\mbox{structure of a positive-dimensional fiber of $f$}&\mbox{singularities of }Y\\
\hline\hline
0&(X,A)\mbox{ is a Fano manifold of co-index five}&\mbox{smooth}\\ \cline{1-3}
1&(F^{n-1} ,A_F )\mbox{ is a Fano manifold of co-index four}&\mbox{smooth}\\ \cline{1-3}
2&(F^{n-2} ,A_F )\mbox{ is a Mukai manifold}&\mbox{$1$-factorial, normal}\\ \cline{1-3}
3&(F^{n-3} ,A_F )\mbox{ is a Del Pezzo manifold}&\mbox{$1$-factorial, normal}\\ \cline{1-3}
4&(F,A_F )\cong(\hyq{n-4},\str(1)),
\mbox{ where }\hyq{n-4}\mbox{ is smooth}
&\mbox{$1$-factorial, normal}\\
\cline{1-3}
5&(F,A_F )\cong(\prs{n-5},\str(1)),\;(n\geq6)&\mbox{$1$-factorial, normal}\\
\cline{1-3}
n&r=0\mbox{ and }(E,A_E )\mbox{ is a Mukai variety such that}&\mbox{$1$-factorial, terminal}\\ 
&\mbox{}\Delta(E,A_E )={A_E}^{n-1} /2\geq1\mbox{ and }\nbl{E/X}\cong-A_E&\\ 
\cline{2-3}
&r=0\mbox{ and }(E,A_E )\mbox{ is a Del Pezzo variety such that}&\mbox{$2$-factorial, terminal}\\ 
&\mbox{}\nbl{E/X}\cong-2A_E&\\ 
\cline{2-3}
&r=0\mbox{ and }(E,A_E ,\nbl{E/X})\cong(\hyq{n-1} ,\str(1),\str(-3))&\mbox{$3$-factorial, terminal}\\
\cline{2-3}
&r=0\mbox{ and }(E,A_E ,\nbl{E/X})\cong(\prs{n-1}, \str(1),\str(-4))&\mbox{$4$-factorial, terminal}\\ 
\cline{2-3}
&r=1\mbox{ and }\rest{\nbl{E/X}}{F}\cong-A_F.
\mbox{ If either $n=5$ or $F$ is irreducible,}
&\mbox{$1$-factorial, terminal}\\
&\mbox{then }(F^{n-2} ,A_F )\mbox{ is a Del Pezzo variety.}&\\
\cline{2-3}
&r=1\mbox{ and }(F,A_F ,\rest{\nbl{E/X}}{F})\cong(\hyq{n-2},\str(1),\str(-2))&\mbox{$2$-factorial, terminal}\\
\cline{2-3}
&r=1\mbox{ and }(F,A_F ,\rest{\nbl{E/X}}{F})\cong(\prs{n-2},\str(1),\str(-3))&\mbox{$3$-factorial, terminal}\\
\cline{2-3}
&r=2\mbox{ and }(F,A_F ,\rest{\nbl{E/X}}{F})\cong(\hyq{n-3},\str(1),\str(-1))&\mbox{$1$-factorial, terminal}\\
\cline{2-3}
&r=2\mbox{ and }(F,A_F ,\rest{\nbl{E/X}}{F})\cong(\prs{n-3},\str(1),\str(-2))&\mbox{$2$-factorial, terminal}\\
\cline{2-3}
&r=3\mbox{ and }(F,A_F ,\rest{\nbl{E/X}}{F})\cong(\prs{n-4},\str(1),\str(-1))&\mbox{$1$-factorial, terminal}\\
\hline
\end{array}$$

\noindent{where} 
$\hyq{}$ 
denotes a hyperquadric in a projective space.
For the definitions of a Fano variety, co-index, a Del Pezzo variety and a Mukai variety, see Definition~\ref{spva}.
\end{TH}

\noindent
If we do not impose the numerical reducedness on the polarization, that is, if ray contractions
$f:X\to Y$
have a supporting divisor
$K_X +(n-4)L$
for a merely ample divisor
$L$,
then it is expected that the list for these ray contractions may contain some extra cases.
Indeed, Proposition~\ref{foray} provides a table for such extra contractions
$(X,L)$
which is recovered from the tables of more ``rigid'' ones, the 
$k$-th adjoint contractions
$(X,A)$
of divisorial type for
$k=1,2$
and
$3$ 
in \cite[Proposition~1.2 and Theorem~1.3]{na:3adc}.

\begin{PR}
\label{foray}
Let
$X$
be an $n$-dimensional smooth projective variety.
Let
$f:X\to Y$
be a ray contraction of divisorial or fiber type.
If a supporting divisor of
$f$
is
$K_X +(n-4)L$
for a merely ample divisor
$L$
on
$X$,
then all the possible cases of
$f$
consist of the following table as well as the table of the fourth adjoint contractions in Theorem~\ref{fmth}.
Throughout the following table,
$f$
is of divisorial type,
$E$
denotes a unique exceptional prime divisor of
$f$,
and
$F$ 
denotes a general fiber of
$\rest{f}{E}:E\to f(E)$
(or any fiber over a smooth point of
$f(E)$
if
$\dim f(E)=1$
and either
$F\cong\prs{4}$
or
$\prs{3}$).

$$\begin{array}{||c|l|c||}\hline
\dim f(E)&\mbox{structure of a positive-dimensional fiber of $f$}&\mbox{singularities of }Y\\ \hline\hline
0&(E,L_E ,\nbl{E/X})\cong(\prs{n-1},\str(9-n),\str(-1))
&\mbox{smooth}\\
&(n=5\mbox{ or }7)&\\
\cline{2-3}
&(E,L_E ,\nbl{E/X})\cong(\prs{n-1},\str(8-n),\str(-2))
&\mbox{$2$-factorial, terminal}\\
&(n=5\mbox{ or }6)&\\
\cline{2-3}
&(E,L_E ,\nbl{E/X})\cong(\hyq{n-1},\str(8-n),\str(-1))
&\mbox{$1$-factorial, terminal}\\
&(n=5\mbox{ or }6)&\\
\cline{2-3}
&(E,L_E ,\nbl{E/X})\cong(\prs{4},\str(2),\str(-3))\;(n=5)&\mbox{$3$-factorial, terminal}\\ 
\cline{2-3}
&(E,L_E ,\nbl{E/X})\cong(\hyq{4},\str(2),\str(-2))\;(n=5)
&\mbox{$2$-factorial, terminal}\\
\cline{2-3}
&(E,L_E ,\nbl{E/X})\cong(E,2A_E ,-A_E )\;(n=5),\mbox{}
&\mbox{$1$-factorial, terminal}\\
&\mbox{where $A$ is a numerically reduced ample divisor}&\\
&\mbox{on $X$, and $(E,A_E )$ is a Del Pezzo variety}&\\ 
\hline
1&(F,L_F ,\rest{\nbl{E/X}}{F})\cong(\prs{n-2},\str(8-n),\str(-1))&\mbox{$1$-factorial, terminal}\\ 
&(n=5\mbox{ or }6)&\\
\cline{2-3}
&(F,L_F ,\rest{\nbl{E/X}}{F})\cong(\prs{3},\str(2),\str(-2))\;(n=5)
&\mbox{$2$-factorial, terminal}\\ 
\cline{2-3}
&(F,L_F ,\rest{\nbl{E/X}}{F})\cong(\hyq{3},\str(2),\str(-1))\;(n=5)
&\mbox{$1$-factorial, terminal}\\ 
\hline
2&(F,L_F ,\rest{\nbl{E/X}}{F})\cong(\prs{2},\str(2),\str(-1))\;(n=5)
&\mbox{$1$-factorial, terminal}\\ 
\hline
\end{array}$$
\end{PR}

\begin{pf}
We follow the Proof of \cite[Proposition~2.1]{na:3adc}.
It's enough to treat only
$f$
of divisorial type.
We can assume that there are a numerically reduced ample divisor
$A$
and an integer
$p\geq2$
such that
$$L\cdot C=pA\cdot C$$
for any curve
$C$
in
$R$.
From [ibid., Lemma~3.3],
$$\dim E+\dim F\geq(n-4)pa+n-1,$$
where
$a=A\cdot C_0$
for some curve
$C_0$
in
$R$
such that
$-K_X \cdot C_0 =l(R)$
(see Definition~\ref{leng} for
$l(R)$).
Since
$\dim E=n-1\geq\dim F$
and
$a\geq1$,
we have
$n-1\geq(n-4)p$
and thus
$$p\leq(n-1)/(n-4)=1+3/(n-4)\leq4.$$
Hence,
$p=2,3$
or
$4$.
If we define an integer
$k$
by
$n-k=(n-4)p$
or
$$k=4p-(p-1)n,$$
then since (after adjusting
$A$
if necessary) we obtain
$R=(K_X +(n-4)pA)^\bot \cap\overline{NE}(X)$,
$f:X\to Y$
can be viewed as the $k$-th adjoint contraction with supporting divisor
$K_X +(n-k)A$.
[ibid., Proposition~1.2] implies that
$k\geq1$
for divisorial type.
Hence,
$n-1\geq(n-4)p$
and thus
$$n\leq(4p-1)/(p-1)=4+3/(p-1).$$

\begin{enumerate}
\item
$p=2$:
$5\leq n\leq7$
and
$k=8-n$.
Hence,
$(k,p,n)=(1,2,7),(2,2,6)$
or
$(3,2,5)$.

\item
$p=3$:
$5\leq n\leq5.5$
and
$k=12-2n=2$.
Hence,
$(k,p,n)=(2,3,5)$.

\item
$p=4$:
$5\leq n\leq5$
and
$k=16-3n=1$.
Hence,
$(k,p,n)=(1,4,5)$.
\end{enumerate}
Therefore, from [ibid., Proposition~1.2 and Theorem~1.3] and
$L_F \cong pA_F$,
we deduce the desired table.
\end{pf}

\noindent
We will obtain Theorem~\ref{frfd} as a corollary of Proposition~\ref{foray}.

\noindent
{\em Proof of Theorem~\ref{frfd}.}
The argument is exactly the same as in the Proof of [ibid., Corollary~2.2] by taking 
$n=5$
in Proposition~\ref{foray}.
\hfill{\em q.e.d.}

\begin{RE}
Suppression of the polarization from the expression of Theorem~\ref{frfd} causes loss of some geometric information, which is visible by comparison with the expression of Proposition~\ref{foray}.
\end{RE}

\noindent
The rest of the paper is devoted to prove Theorem~\ref{fmth}.
In Sections $3$ through $7$, we apply the techniques developed in \cite{na:3adc}, in particular, the Generality on the 
$k$-th adjoint contractions in [ibid., Section~4] in order to determine a general positive-dimensional fiber
$F$
of the fourth adjoint contraction
$f:X\to Y$
of divisorial and fiber types.
However, it is far beyond these techniques of [ibid.] to show that a general positive-dimensional fiber 
$F$
is irreducible for
$f:X\to Y$
of divisorial type with
$\dim f(E)=1$
and
$\rest{\nbl{E/X}}{F}\cong-A_F$,
where
$E$
is a unique exceptional prime divisor of
$f$.
If we knew that
$E$
is normal, it would be easy to show that a general fiber
$F$
of
$\rest{f}{E}:E\to f(E)$,
the restriction of
$f$
to
$E$,
is normal and so irreducible since
$f$ 
has only connected fibers.
Unfortunately, we do not know whether
$E$
is, in general, normal, and this is an open question.

\begin{OQ}
\label{epdn}
Is a unique exceptional prime divisor of a ray contraction of divisorial type always normal?
\end{OQ}

\noindent
Aside from the context of ray contraction, we can construct the following example which might ``tempt'' a possibility of counterexample against the affirmative answer to this Question~\ref{epdn}.

\begin{EX}
(Z. Ran)
Let
$\pi:\prs{1}\times\prs{1}\to\prs{1}$
be a natural projection and
$\rho:\prs{1}\to\prs{1}$
a double cover defined by
$x\to x^2$
in a local coordinate of
$\prs{1}$.
Let
$C$
be the image of a section of
$\pi$.
Note that
$C$
intersects at exactly two points (counting the multiplicity) with each fiber of
$h:=\rho\circ\pi$.
By identifying these two points on
$C$,
we construct a variety
$E$
and a morphism
$g:E\to\prs{1}$
from
$h:\prs{1}\times\prs{1}\to\prs{1}$.
Then,
$E$
is irreducible and a general fiber of
$g$
is connected and reducible.
\end{EX}

\noindent
In \cite[Proof of Theorem~2.1]{an:erhdv}, T. Ando asserts that Bertini Theorem would prove that 
$F$
is irreducible,
but the author has not figured out how to apply Bertini Theorem to the possibly non-normal variety
$E$
and this does not seem a trivial application of this theorem.

\smallskip

\noindent
The method we will use in Section~\ref{pirr} is totally different from Ando's and does not depend on Bertini Theorem as explained in what follows:
Let
$\rest{f}{E}:E\to f(E)$
and
$F$
be as above.
Let
$M$
be the pull back by
$f$
of a general hyperplane section
$V$
on
$Y$.
Then,
$F$
is realized as a connected component of the exceptional locus of
$\rest{f}{M}:M\to V$.
Assuming
$F$
is not irreducible, we will eventually derive a contradiction when
$n:=\dim X=5$.
Let
$G$
be an arbitrary irreducible component of 
$F$,
and
$\nu:=\dim G \; (\nu=n-2=3)$.
We understand
$F$
and
$G$
as divisors on
$M$.
Then, from the assumption there is ``an information divisor''
$D\; (>0)$
of
$G$
which is the restriction to
$G$
of the other components of 
$F$
and tells how 
$G$
intersects with the other components.
We will list up
$(\tilde{G},A_{\tilde{G}} ,\nbl{G/M},D)$
completely in Lemma~\ref{clofgand}, where
$(\tilde{G},A_{\tilde{G}} )$
is the normalization of
$(G,A_G )$.
Then, we will calculate the intersection number
$F\cdot Z$
with a fixed curve
$Z$
in two different methods:
The first one,
$F\cdot Z=\sum_G G\cdot Z$
based on this list.
Secondly,
$F\cdot Z$
is calculated from the exceptional divisor
$E$
of 
$f$.
As a result, these two methods produce two distinct numbers for the same
$F\cdot Z$,
hence, the expected contradiction.
(See the beginning in Section~\ref{pirr} for the precise notations.)

\smallskip

\noindent
The most difficult part to obtain is the list of
$(\tilde{G},A_{\tilde{G}} )$
when
$\dim f(E)=1$
and
$\rest{\nbl{E/X}}{F}\cong-A_F$
(Lemma~\ref{clofg}).
In order to apply T. Fujita's classification theory, we have to show
$\Delta(\tilde{G},A_{\tilde{G}} )=0$
in Lemma~\ref{fstsdg}.
(See Lemma~\ref{nong} for
$\Delta$-genus.)
The key tool to prove this equation is the sectional genera (Definition~\ref{sege}),
$g(G,A_G )$
and
$g(\tilde{G},A_{\tilde{G}} )$
which behave well with singularities.
We will actually generalize the sectional genus formula, as evidence of its good behavior, to Cohen-Macaulay pre-polarized varieties in Proposition~\ref{sgf}.
From this Proposition~\ref{sgf}, the condition
$D>0$
implies that
$g(G,A_G )\leq0$
(Lemma~\ref{fstsdg}).
Then, another useful Proposition~\ref{inse}, measurement of non-normal locus by the sectional genus, asserts that
$$g(G,A_G )\geq g(\tilde{G},A_{\tilde{G}} ).$$
Therefore, we obtain
$g(\tilde{G},A_{\tilde{G}} )\leq0$.
What we need here is exactly Fujita's conjecture regarding sectional genus (\cite[Conjecture in Introduction]{fj:rqp}).

\begin{CON}
\label{fcsg}
It should be true that
$g\geq0$
for any quasi-polarized projective variety (Definition~\ref{pv}).
Moreover,
$g=0$
should imply
$\Delta=0$
if the variety is normal.
\end{CON}

\noindent
Only for normal varieties of dimension three or less, this conjecture is solved by Fujita in \cite[Corollary~(4.8)]{fj:rqp}, where he essentially reduces the non-negativity of the sectional genus to the nefness of canonical divisor, or the existence of minimal models in the $3$-dimensional Minimal Model Program in \cite{mo:3mmp}.
Therefore, we are able to obtain the list of
$(\tilde{G},A_{\tilde{G}} )$
only when
$\nu=\dim \tilde{G}=3$
or equivalently
$n=\dim X=5$,
while the unsolved part of Fujita's Conjecture~\ref{fcsg} remains in our setting as the following

\begin{CON}
A general positive-dimensional fiber
$F$
of the fourth adjoint contraction of divisorial type should be irreducible when
$\dim f(E)=1$
and
$\rest{\nbl{E/X}}{F}\cong-A_F$
even for
$\dim X>5$.
\end{CON}

\noindent
The construction of each case of the lists is technically rather different business.
We just leave the following

\begin{PROB}
Construct each case of the lists in Theorems~\ref{frfd}, \ref{fmth} and Proposition~\ref{foray} if it exists.
\end{PROB}

\noindent
It should be mentioned that the following case be included in \cite[Proposition~1.2]{na:3adc}.

$$\begin{array}{||c|c|c|c||}\hline
k&\dim Y&
\mbox{structure of a positive-dimensional fiber of $f$}
&\hspace{.1in}\mbox{singularities of }Y\hspace{.1in}\\ \hline\hline
\vspace*{\fill}1\vspace*{\fill}&2&\mbox{$f$ is a scroll $(n\geq3)$}&\mbox{smooth}\\ \cline{2-4} \hline
\end{array}$$

\medskip

\noindent{\em Acknowledgments.}
A part of this article is written when the author was a Visiting Lecturer at Oklahoma State University, August 1994 through May 1995.
Thanks are due to Professor Sheldon Katz for the invitation to research there.
The completion for the five dimensional case was done when the author was a Visiting Assistant Professor at the University of California, Riverside, September 1997 through June 1998.
Thanks are also due to Professor Mei-Chu Chang and in particular Professor Ziv Ran for many helpful discussions.
The main result, Theorem~\ref{frfd} is presented along with the outline of its proof in {\em A Conference on Algebraic Geometry to celebrate Robin Hartshorne's 60th birthday} at UC Berkeley on August 30, 1998.
The author would like to thank the organizers for the support.

\medskip

\section{Basic concepts}
\label{baco}

\begin{DE}
A {\em variety} means a separated integral scheme of finite type over the complex
number field 
$\Bbb C$
unless stated otherwise. 
A {\em subvariety} means an integral subscheme of a scheme. 
A {\em point} means a closed point.
Locally free sheaves and vector bundles are used interchangeably.
$\dsh{X}$
denotes the dualizing sheaf of a complete scheme
$X$.
$K_X$
denotes the canonical divisor or the canonical sheaf of a normal variety
$X$.
A normal variety is said to have $q$-{\em factorial singularities} 
if there is an integer
$q\geq1$
such that for every Weil divisor 
$Z$,
$qZ$ 
is a Cartier divisor.
For convenience, $0$-{\em factorial} means non-singular.
(In our case, $1$-factoriality is equivalent to local factoriality, as shown in 
\cite[Proposition~(3.10) on p.\ 139]{ak:duality}.)
We employ the intersection theory due to Kleiman in \cite{kl:amp} as well as Definition~\ref{adfm}.
The linear equivalence between divisors
$D_1$
and
$D_2$
is denoted by
$D_1 \sim D_2$.
\end{DE}

\begin{DE}
\label{leng}
Let 
$X$ 
be an 
$n$-dimensional 
smooth projective variety with an extremal ray 
$R$. 
Let 
$f: X\to Y$ 
be the ray contraction of 
$R$ 
as in \cite[Definition~3-2-3]{kmm}.
Set
$$\Etil=\{x\in X\mid f\mbox{
is not isomorphic at 
}x\}.$$
Then, we call the ray contraction 
$f$ 
to be {\em of fiber type}, {\em of divisorial type}, or {\em of flipping type} 
(equivalently, {\em small}) if 
$\dim\Etil=n$, 
$=n-1$, 
or 
$\leq n-2$, 
respectively, 
where 
`$\dim $' 
denotes the maximum of the dimensions of the irreducible components of
$\Etil$.
We define the {\em length} of the extremal ray
$R$
by
$$l(R)=\min\{-K_X \cdot C; C\mbox{ 
is a possibly singular rational curve whose numerical class belongs to 
}R\}.$$
\end{DE}

\smallskip

\noindent
We will make small contributions, two Propositions~\ref{inse} and \ref{sgf}, to the Theory of Polarized Varieties after T. Fujita.

\begin{DE}(cf.\  \cite[Introduction]{fj:rqp})
\label{pv}
Let
$X$
be a complete variety over an algebraically closed field, and
$L$
a line bundle on it.
A pair
$(X,L)$
is called
a {\em polarized (quasi-polarized, pre-polarized) variety} if
$L$
is an ample (a nef and big, a merely) line bundle, respectively.
\end{DE}

\begin{LE}([ibid., Theorem~(1.1)])
\label{nong}
Let
$\Delta(X,L)$
be the {\em delta genus} of an $n$-dimensional quasi-polarized variety
$(X,L)$, i.e.
$$\Delta(X,L)=n+L^n -h^0 (X,L).$$
If
$X$
is a projective variety, then
$\Delta(X,L)$
is a non-negative integer.
\end{LE}

\begin{DE}(cf.\  [ibid., Introduction])
\label{sege}
Let
$(X,L)$
be an $n$-dimensional pre-polarized variety.
We define integers
$\chi_i$
by
$$\chi(X,tL)=\sum_{i=0}^{n} (\chi_i /i!)t^{[i]} ,$$
where
$t^{[i]}=t(t+1)\cdots(t+i-1)$
for
$i>0$
and
$t^{[0]} =1$.
Then, the {\em sectional genus}
$g(X,L)$
of
$(X,L)$
is defined by
$$g(X,L)=1-\chi_{n-1} \;(\in\integer).$$
\end{DE}

\noindent
The sectional genus detects if the size of the non-normal locus of a polarized variety is large.
We will need the following proposition to prove Lemma~\ref{fstsdg}.

\begin{PR}
\label{inse}
Let
$\mu:\tilde{X}\to X$
be the normalization of a complete variety
$X$
over an algebraically closed field and
$L$
an ample line bundle on
$X$.
Then, we have
\begin{enumerate}
\item
$g(\tilde{X},\mu^* L)=g(X,L)$
if the non-normal locus of
$X$
has codimension
$>1$, and
\item
$g(\tilde{X},\mu^* L)<g(X,L)$
if the non-normal locus of
$X$
has codimension
$=1$.
\end{enumerate}
\end{PR}

\begin{pf}
We will go through a Leray spectral sequence argument.
Let
$n=\dim X$.
Let
$E^{p,q}_2 =\Coh{p}{X}{R^q\mu_{\ast} \mu^* (tL)}$.
Since
$R^q\mu_{\ast} \mu^* (tL)=0$
for
$q>0$,
we have
\begin{equation*}
E^{p,q}_2 =
\begin{cases}
\Coh{p}{X}{\mu_{\ast} \str_{\tilde{X}} \otimes tL} & \text{if $q=0$},\\
0&\text{if $q>0$}.
\end{cases}
\end{equation*}
Hence,
$E^{p,q}_2 \cong E^{p,q}_{\infty}$
and thus
$$\Coh{i}{\tilde{X}}{\mu^* (tL)}\cong\oplus_{p+q=i} E^{p,q}_{\infty}
=E^{i,0}_{\infty} =\Coh{i}{X}{\mu_{\ast} \str_{\tilde{X}} \otimes tL}.$$
Therefore, we obtain the following relation between two Euler characteristics on
$X$
and
$\tilde{X}$.
$$\chi(\tilde{X},t\mu^* L)=\chi(X,\mu_{\ast} \str_{\tilde{X}} \otimes tL).\eqno(\ref{baco}.1)$$
There is a natural short exact sequence,
$\shortseq{\str_X}{}{\mu_{\ast} \str_{\tilde{X}}}{}{\cal F}$,
where
${\cal F}$
is a coherent torsion sheaf.
Hence, form (\ref{baco}.1),
$$\chi(\tilde{X},t\mu^* L)=\chi(X,tL)+\chi(X,{\cal F}\otimes tL).\eqno(\ref{baco}.2)$$
Since
$\dim\supp{\cal F}\leq n-1$,
$$\chi(X,{\cal F}\otimes tL)=[a_{n-1} /(n-1)!]t^{n-1} +\cdots+a_0 \quad(a_i \in\rational).\eqno(\ref{baco}.3)$$
From Definition~\ref{sege}, (\ref{baco}.2) and (\ref{baco}.3), we deduce
$g(\tilde{X},\mu^* L)=g(X,L)-a_{n-1}$.
Therefore, to finish the proof, it suffices to show the following Claim~\ref{apoz} since the non-normal locus coincides with
$\supp{\cal F}$.
\begin{CL}
\label{apoz}

\begin{equation*}
a_{n-1}
\begin{cases}
=0&\text{if $\dim\supp{\cal F}<n-1$},\\
>0&\text{if $\dim\supp{\cal F}=n-1$}.
\end{cases}
\end{equation*}
\end{CL}

\noindent
This is a direct consequence from
\cite[Proposition~1 on p.\ 302]{kl:amp}.
Since
$L$
is ample, putting
$s=\dim\supp{\cal F}$,
this Proposition~1 implies
$$\chi(X,{\cal F}\otimes tL)=[\deg{\cal F}/s!]t^s +(\mbox{terms of lower degree}),$$
and
$\deg{\cal F}>0$.
Hence, from (\ref{baco}.3), it follows that
$a_{n-1} =0\mbox{ or }\deg{\cal F}$
if
$s<n-1\mbox{ or }=n-1$,
respectively.
\end{pf}

\noindent
The following definition provides the intersection number of a coherent sheaf of rank one with an invertible sheaf.
Note that the projectiveness imposed on a variety in \cite[Definition~3.6]{na:3adc} is loosened to the completeness, and the amplitude on an invertible sheaf in [ibid.] is dropped in the new definition.

\begin{LEDE}
\label{adfm}
Let
$X$
be a complete variety over an algebraically closed field, of dimension
$n\geq2$.
Let
$\cal F$
be any coherent sheaf of rank one, and let
$A$
and
$B$
be any invertible sheaves on
$X$.

\begin{enumerate}
\item
Then, we obtain
$$\euch{X}{{\cal F}\otimes B^{\otimes t}}
=(B^n /n!)t^n +O(t^{n-1} ).$$
Therefore, we can define the {\em intersection number}
${\cal F}\cdot B^{n-1}$,
$\in\integer$,
of a coherent sheaf
${\cal F}$
of rank one with an invertible sheaf
$B$
by the following equation
$$\euch{X}{{\cal F}\otimes B^{\otimes t}}-\euch{X}{B^{\otimes t}}
=({\cal F}\cdot B^{n-1} /(n-1)!)t^{n-1} +O(t^{n-2} ).$$

\item
Moreover, applying this definition of intersection number, we recover partially \cite[Proposition~2 on p.\ 296]{kl:amp} in the following form
$$({\cal F}\otimes A)\cdot B^{n-1}
={\cal F}\cdot B^{n-1} +A\cdot B^{n-1} .$$ 
\end{enumerate}
\end{LEDE}

\begin{pf}
Note that
$\dim \supp {\cal F}=n$
for a coherent sheaf of rank one.
From \cite[Theorem on p.\ 295, Lemma~1 on p.\ 296]{kl:amp}, there are integers
$a_{i,j}$
and
$b_i$
such that for any integers
$s$
and
$t$,
$$\begin{array}{rcl}
\euch{X}{{\cal F}\otimes A^{\otimes s} \otimes B^{\otimes t}}
&=&
\sum_{i,j\geq0, \; i+j\leq n} a_{i,j}
\left(
\begin{array}{c}
s+i \\
i
\end{array}
\right)
\left(
\begin{array}{c}
t+j \\
j
\end{array}
\right) \\
&=&
(a_{0,n} /n!)t^n +[1/(n-1)!][(n+1)a_{0,n} /2+a_{0,n-1} +a_{1,n-1} s]t^{n-1} +O(t^{n-2} )\\
&&\hfill(\ref{baco}.a)\\
\euch{X}{B^{\otimes t}}
&=&
\sum_{i=0}^n b_i
\left(
\begin{array}{c}
t+i \\
i
\end{array}
\right) \\
&=&
(b_n /n!)t^n +[1/(n-1)!][b_{n-1} +(n+1)b_n /2]t^{n-1} +O(t^{n-2} ) \hfill(\ref{baco}.b)
\end{array}
$$
Take
$s=0$
in 
$(\ref{baco}.a)$,
and
$$\euch{X}{{\cal F}\otimes B^{\otimes t}}
=
(a_{0,n} /n!)t^n +[1/(n-1)!][a_{0,n-1} +(n+1)a_{0,n} /2]t^{n-1} +O(t^{n-2} ) \hfill\eqno(\ref{baco}.c)
$$
Substituting
$t=t_1 +\cdots+t_{n-1}$
for 
$(\ref{baco}.a),$
and
$t=t_1 +\cdots+t_n$
for
$(\ref{baco}.b)$ 
and
$(\ref{baco}.c),$
Kleiman's definition of intersection number in \cite[p.\ 296]{kl:amp} implies that
$$a_{1,n-1} =A\cdot B^{n-1} , b_n =B^n  \mbox{ and } a_{0,n} =(B^n \cdot {\cal F}),\eqno(\ref{baco}.d)$$
where
$(B^n \cdot {\cal F})$
is Kleiman's notation, and this will be shown to be
$B^n$
as follows:
Since
$X$
is reduced and irreducible, and
$\cal F$
is a coherent sheaf of rank one,
$$\mbox{length}_{\str_{X,\xi}} {\cal F}_\xi =1,$$
where
$\xi$
is the generic point of
$X$.
Hence, from [ibid., Corollary~2 on p.\ 298] and
$(\ref{baco}.d),$
$$\begin{array}{rcl}
\quad\quad\quad\quad\quad\quad\quad\quad\quad\quad\quad\quad\quad\quad
a_{0,n} 
&=&(B^n \cdot {\cal F}) \\
&=&(\mbox{length}_{\str_{X,\xi}} {\cal F}_\xi )(B^n \cdot\str_X ) \quad\quad\quad\quad\quad\quad\quad\quad\quad\quad\quad\quad\quad\quad \\
&=&B^n . \hfill (\ref{baco}.e)
\end{array}
$$
Therefore, we obtain the first statement from
$(\ref{baco}.c).$

\noindent
Secondly, take
$s=1$
in 
$(\ref{baco}.a),$ 
and
$$
\euch{X}{{\cal F}\otimes A\otimes B^{\otimes t}}
=
(a_{0,n} /n!)t^n +[1/(n-1)!][a_{0,n-1} +a_{1,n-1} +(n+1)a_{0,n} /2]t^{n-1} +O(t^{n-2} ) \eqno(\ref{baco}.f)
$$
From 
$(\ref{baco}.d)$ 
and 
$(\ref{baco}.e),$
$b_n =a_{0,n}$.
Thus, 
$(\ref{baco}.b),$
$(\ref{baco}.c)$
and 
$(\ref{baco}.f)$ 
imply that
$$\begin{array}{rcl}
\euch{X}{{\cal F}\otimes A\otimes B^{\otimes t}}-\euch{X}{B^{\otimes t}}
&=&
[1/(n-1)!](a_{0,n-1} +a_{1,n-1} -b_{n-1} )t^{n-1} +O(t^{n-2} ),
\end{array}$$
and
$$\begin{array}{rcl}
\euch{X}{{\cal F}\otimes B^{\otimes t}}-\euch{X}{B^{\otimes t}}
&=&
[1/(n-1)!](a_{0,n-1} -b_{n-1} )t^{n-1} +O(t^{n-2} ).
\end{array}$$
Therefore, applying the above Definition~\ref{adfm}, we have
$$\begin{array}{rcl}
\quad\quad\quad\quad\quad\quad\quad\quad\quad\quad
({\cal F}\otimes A)\cdot B^{n-1}
&=&
(a_{0,n-1} -b_{n-1} )+a_{1,n-1} ,
\quad\quad\quad\quad\quad\quad\quad\quad\quad\quad\quad\quad\quad\\
& & \hfill (\ref{baco}.g) \\
{\cal F}\cdot B^{n-1}
&=&
a_{0,n-1} -b_{n-1} .
\end{array}$$
Hence,
$(\ref{baco}.d)$ 
and 
$(\ref{baco}.g)$ 
will provide the second statement.
Also,
$(\ref{baco}.g)$ 
shows that
${\cal F}\cdot B^{n-1}$
is an integer.
\end{pf}

\smallskip

\noindent
We will generalize the sectional genus formula to a Cohen-Macaulay pre-polarized variety and use this in the Proofs of Lemmas~\ref{fstsdg} and \ref{gnml}.
 
\begin{PR}
\label{sgf}
Let
$L$
be an invertible sheaf on an $n$-dimensional complete variety
$X$
over an algebraically closed field
$(n\geq2)$.
If
$X$
is locally Cohen-Macaulay, then the sectional genus formula holds:
$$2g(X,L)-2=[\dsh{X}\otimes L^{\otimes(n-1)}]\cdot L^{n-1} ,$$
where we use the intersection number defined in Definition~\ref{adfm}.
In particular, 
if some additional conditions allow to express
$\dsh{X}\otimes L^{\otimes(n-1)}$
additively as a Weil divisor, 
for instance, 
either
$\dsh{X}$
is invertible or
$\dsh{X}$
is a reflexive sheaf on a normal variety
$X$,
then we obtain the usual form of the sectional genus formula:
$$2g(X,L)-2=[\dsh{X}+(n-1)L]\cdot L^{n-1} .$$
\end{PR}

\begin{pf}
Since
$X$
is proper over an algebraically closed field, the dualizing sheaf
$\dsh{X}$
exists (see \cite{gr:abv} and \cite{ha:redu}).
From Definition~\ref{sege},
$$\begin{array}{rcl}
\quad\quad\quad
\euch{X}{L^{\otimes t}}
&=&
\sum_{i=0}^n (1/i!){\chi_i}\: t^{[i]} \\
&=&
(1/n!){\chi_n}\: t^n +[1/2(n-1)!][(n-1)\chi_n +2\chi_{n-1} ]t^{n-1} +O(t^{n-2}) 
\quad\quad\quad\quad(\ref{baco}.4)
\end{array}$$
From the Serre Duality Theorem for complete varieties (see [ibid.]),
$$\begin{array}{rcl}
\euch{X}{\dsh{X}\otimes L^{\otimes t}}
&=&
(-1)^n \euch{X}{L^{\otimes(-t)}} \\
&=&
(1/n!)\chi_n t^n -[1/2(n-1)!][(n-1)\chi_n +2\chi_{n-1} ]t^{n-1} +O(t^{n-2})
\end{array}$$
Hence,
$$\euch{X}{\dsh{X}\otimes L^{\otimes t}}-\euch{X}{L^{\otimes t}}=-\frac{(n-1)\chi_n +2\chi_{n-1}}{(n-1)!} t^{n-1} +O(t^{n-2}). \eqno(\ref{baco}.5)$$
Therefore, from Definition~\ref{sege} and Definition~\ref{adfm},
$$\begin{array}{rcl}
\dsh{X}\cdot L^{n-1}&=&-(n-1)\chi_n -2\chi_{n-1} \\
                    &=&-(n-1)L^n -2(1-g(X,L)).
\end{array}$$
Hence,
$$2g(X,L)-2=\dsh{X}\cdot L^{n-1}+(n-1)L^n ,$$
and the second statement in Lemma~\ref{adfm} will complete the proof by substituting
${\cal F}=\dsh{X}$
and
$A=B=L$.
\end{pf}

\smallskip

\noindent
The following is the key lemma to prove the irreducibility of a general positive-dimensional fiber of the fourth adjoint contraction of divisorial type when its dimension is five.

\begin{LE}\cite[Corollary~(4.8)]{fj:rqp}
\label{pssg}
Let 
$(V,L)$
be a quasi-polarized variety over an algebraically closed field of characteristic zero.
If
$V$ 
is normal and
$\dim V\leq3$,
then its sectional genus satisfies
$g(V,L)\geq0$.
Moreover, if 
$g(V,L)=0$,
it follows that
$\Delta(V,L)=0$.
\end{LE}

\smallskip

\noindent
Finally, we will clarify the definitions of special varieties arising in this paper.

\begin{DE}
\label{spva}
A projective variety
$X$
with Gorenstein singularities is said to be a {\em Fano variety} if
$-\dsh{X}$
is ample.
If
$-\dsh{X}\sim mL$
(linear equivalence) for an ample Cartier divisor
$L$
on
$X$
and an integer
$m>0$,
then we call
$c=n+1-m$
the {\em co-index} of a Fano variety
$X$.
In particular, a Fano variety of co-index two or three is called a {\em Del Pezzo variety} or a {\em Mukai variety}, respectively.
\end{DE}

\begin{RE}
\label{nodf}
In the above Definition~\ref{spva}, for simplicity, we have relaxed the definition of co-index which is defined in \cite[Definition~3.8]{na:3adc}.
Actually, among all the lists in this paper, only a Fano manifold
$(X,A)$ 
of co-index five in Theorem~\ref{fmth} is known to be of co-index five in the strict sense of [ibid.].
\end{RE}

\begin{LE}
\label{mdfn}
Let
$X$
be a Fano variety which has log-terminal singularities (over algebraically closed field of characteristic zero).
If 
$-\dsh{X}\sim mL$
for an ample Cartier divisor
$L$
on
$X$
and an integer
$m>0$,
then,
we have
$$\Coh{i}{X}{tL}=0\mbox{ for any $t\in\integer$ and any $0<i<\dim X$}.
\eqno(\ref{baco}.6)$$
\end{LE}

\begin{pf}
From the Serre Duality Theorem, it follows that
$\Coh{i}{X}{tL}=\Coh{n-i}{X}{\dsh{X}-tL}\;(n=\dim X)$.
Apply the Kawamata-Viehweg Vanishing Theorem, \cite[Theorem~1-2-5 and Remark~1-2-6]{kmm} to each side of this equality.
\end{pf}

\begin{RE}
\label{cvgf}
Let
$f:X\to Y$
be an $n$-dimensional fourth adjoint contraction of divisorial type and 
$F$
a positive general fiber of
$f$.
From \cite[Section~4]{na:3adc}, we deduce that 
$$\Coh{i}{F}{tA_F}=0
\mbox{ for any
$0<i<\dim X$
and any $t\in\integer$ satisfying either
$t\leq-q$
or
$-(n-4)\leq t$}$$
(see [ibid.] for the notations).
We will see later that
$q$
takes only
$1$
or
$2$
when Del Pezzo and Mukai varieties arise as
$F$.
(They appear only in Sections~\ref{prfz} and \ref{prfo}.)
Therefore, together with the above Lemma~\ref{mdfn}, we will conclude that all the Del Pezzo and Mukai varieties in every list in this paper satisfy (\ref{baco}.6), or the traditional definitions of these varieties as in \cite[Definition~3.8]{na:3adc} and \cite[Section~2]{an:erhdv}.
\end{RE}

\medskip

\section{Numerical classification of the fourth adjoint contraction of divisorial type and the proof of fiber type}
\label{nufi}

\begin{PR}
\label{nucl}
If
$f:X\to Y$ 
is the fourth adjoint contraction of divisorial type, then we have
$r=\dim f(E)=0,1,2$
or
$3$.
\end{PR}

\begin{pf}
Exactly the same argument as in the Proof of \cite[Proposition~6.1]{na:3adc} works for this case.
\end{pf}

\begin{pf*}{Proof of fiber type}
Let
$F$
be a general fiber of
$f:X\to Y$.
Then, applying a Bertini Theorem, \cite[Corollary~10.9 and Remark~10.9.1 in Chap.\ III]{ha:alg}, it follows that
$F$
is smooth and irreducible since
$f$
has only connected fibers.
As in \cite[Section~6]{na:3adc}, we have
$\dim F=n,n-1,n-2,n-3,n-4,n-5$
and
$K_F +(n-4)A_F \cong\str_F$.
Hence, Definition~\ref{spva} will determine
$F$.
The Proof of \cite[Lemma~5-1-5]{kmm} works for our purpose and produces the statement of singularities of
$Y$.
\end{pf*}

\medskip

\section{Proof of the fourth adjoint contraction of divisorial type with
$\protect{\dim}f(E)=0$}
\label{prfz}

\noindent{In} this case,
$(F,A_F )$
is a polarized variety
$(E,A_E )$.
We will use the same argument as in \cite[Section~7]{na:3adc}.
Let
$\nu=\dim E\;(=n-1\geq4)$,
$d={A_E}^\nu$
and
$P(t)=\chi(E,tA_E )$.
From [ibid., Lemma~4.5], it follows that
$$P(t)=(d/\nu!)(t+1)\cdots(t+\nu-3)(t+\alpha)(t+\beta)(t+\gamma)\eqno(\ref{prfz}.1)$$
for some complex numbers
$\alpha$,
$\beta$
and
$\gamma$,
$$\nu(\nu-1)(\nu-2)=d\alpha\beta\gamma,\eqno(\ref{prfz}.2)$$
$$\coh{0}{E}{A_E}=[d/\nu(\nu-1)][1+(\alpha+\beta+\gamma)+(\beta\gamma+\gamma\alpha+\alpha\beta)+\alpha\beta\gamma],\eqno(\ref{prfz}.3)$$
and
$$P(t)=(-1)^\nu P(-t-(\nu-3)-q).\eqno(\ref{prfz}.4)$$
From (\ref{prfz}.1) and (\ref{prfz}.4), we have two expressions for roots of
$P(t)=0$:
$$\{-1\geq\cdots\geq-(\nu-3),-\alpha,-\beta,-\gamma\}=$$
$$\{-1+(1-q)\geq\cdots\geq-(\nu-3)+(1-q),\alpha-(\nu-3)-q,\beta-(\nu-3)-q,\gamma-(\nu-3)-q\}.\eqno(\ref{prfz}.5)$$
By adding the both sides of (\ref{prfz}.5),
$$\alpha+\beta+\gamma=\nu-3+q\nu/2.\eqno(\ref{prfz}.6)$$
By multiplying the both sides of (\ref{prfz}.5),
$$(-1)^\nu (\nu-3)!\alpha\beta\gamma=(-1)^{\nu-3} [(\nu+q-4)!/(q-1)!]$$
$$\cdot[\alpha\beta\gamma-(\nu-3+q)(\beta\gamma+\gamma\alpha+\alpha\beta)+(\nu-3+q)^2 (\alpha+\beta+\gamma)-(\nu-3+q)^3 ].$$
Hence,
$$\beta\gamma+\gamma\alpha+\alpha\beta=[1/(\nu-3+q)][1+(q-1)!(\nu-3)!/(\nu+q-4)!]\alpha\beta\gamma$$
$$+(\nu-3+q)(\alpha+\beta+\gamma)-(\nu-3+q)^2 .\eqno(\ref{prfz}.7)$$
Substitute (\ref{prfz}.2), (\ref{prfz}.6) and (\ref{prfz}.7) into (\ref{prfz}.3), and then
$$\coh{0}{E}{A_E}=[d/\nu(\nu-1)]\{1-(\nu-3+q)^2 +(\nu-2+q)(\nu-3+q\nu/2)$$
$$+[1/(\nu-3+q)][\nu-2+q+(q-1)!(\nu-3)!/(\nu+q-4)!][\nu(\nu-1)(\nu-2)/d]\}$$
$$=[(\nu-2+q)(q\nu-2q+2)/(2\nu(\nu-1))]d$$
$$+[(\nu-2)/(\nu-3+q)][\nu-2+q+(q-1)!(\nu-3)!/(\nu+q-4)!].\eqno(\ref{prfz}.8)$$
[ibid., Lemma~4.3] implies
$$\dsh{E}\cong-(\nu-3+q)A_E .\eqno(\ref{prfz}.9)$$
Hence, from [ibid., Lemma~5.1], we have
$\nu-3+q\leq\nu+1$,
and thus
$q\leq4$.

\noindent{We} will examine each of four cases:
$q=1,2,3,4$.

\medskip

{\bf Case 1)}
$q=1$:
From \cite[Lemma~4.6]{na:3adc},
$Y$
has $1$-factorial, terminal singularities, and [ibid., Lemma~4.2] implies that
$\nbl{E/X}\cong-A_E$.
From (\ref{prfz}.8), we deduce
$\coh{0}{E}{A_E}=d/2+\nu$
and thus
$\Delta(E,A_E )=d/2$.
From (\ref{prfz}.9),
$\dsh{E}\cong-(\nu-2)A_E$.
Since
$E$
has Gorenstein singularities,
$(E,A_E )$
is a Mukai variety such that
$\Delta(E,A_E )=d/2\geq1$.
(See Definition~\ref{spva} and Remark~\ref{nodf} for the relation with the traditional definition of Mukai variety.)

\medskip

{\bf Case 2)}
$q=2$:
As above
$Y$
has
$2$-factorial, terminal singularities, and
$\nbl{E/X}\cong-2A_E$.
From (\ref{prfz}.9),
$\dsh{E}\cong-(\nu-1)A_E$.
Since
$E$
has Gorenstein singularities,
$(E,A_E )$
is a Del Pezzo variety.
(See Definition~\ref{spva} and Remark~\ref{nodf} for the relation with the traditional definition of Del Pezzo variety.)

\medskip

{\bf Case 3)}
$q=3$:
$Y$ has $3$-factorial, terminal singularities, and
$\nbl{E/X}\cong-3A_E$.
From (\ref{prfz}.8),
$$\coh{0}{E}{A_E}=[(3\nu^2 -\nu-4)d+(2\nu^3 -4\nu^2 -2\nu+8)]/2\nu(\nu-1)$$
and thus
$$\Delta(E,A_E )=(2-d)(\nu^2 -\nu-4)/2\nu(\nu-1).\eqno(\ref{prfz}.10)$$
Since
$\nu\geq4$,
the non-negativity of delta genus in Lemma~\ref{nong} implies that
$d=1$
or
$2$.
If
$d=1$,
then from (\ref{prfz}.10),
$\nu[(\nu+1)-2(\nu-1)\Delta(E,A_E )]=4$.
From this, we would have
$\nu=4$
and
$3\Delta(E,A_E )=2$,
which contradicts that a delta genus is an integer.
Hence,
$d=2$
and
$\Delta(E,A_E )=0$.
Therefore, from \cite[Corollary~(4.3)]{fj:sdt}, it follows that
$(E,A_E ,\nbl{E/X})\cong(\hyq{n-1},\str(1),\str(-3))$.

\medskip

{\bf Case 4)}
$q=4$:
$Y$ has $4$-factorial, terminal singularities, and
$\nbl{E/X}\cong-4A_E$.
From (\ref{prfz}.8),
$$\coh{0}{E}{A_E}=[(2\nu^2 +\nu-6)d+(\nu^3 -2\nu^2 -2\nu+6)]/\nu(\nu-1)$$
and thus
$$\Delta(E,A_E )=(1-d)(\nu^2 +2\nu-6)/\nu(\nu-1).$$
As above
$d=1$
and
$\Delta(E,A_E )=0$.
Therefore, from \cite[Corollary~(4.3)]{fj:sdt}, we have that
$(E,A_E ,\nbl{E/X})\cong(\prs{n-1},\str(1),\str(-4))$.

\medskip

\section{Proof of the fourth adjoint contraction of divisorial type with
$\protect{\dim}f(E)=1$}
\label{prfo}

\noindent 
We begin with applying the same argument as in \cite[Section~8]{na:3adc}.
Let
$F$
be a general fiber of
$E\to f(E)$,
$\nu=\dim F \;(=n-2\geq3),\;d={A_F}^\nu$
and
$P(t)=\chi(F,tA_F )$.
From \cite[Lemma~4.5]{na:3adc} we have
$$P(t)=(d/\nu!)(t+1)\cdots(t+\nu-2)(t+\alpha)(t+\beta)$$
for some complex numbers
$\alpha$
and
$\beta$,
$$\nu(\nu-1)=d\alpha\beta,$$
$$\coh{0}{F}{A_F}=d/\nu+(d/\nu)(\alpha+\beta)+(\nu-1),$$
and
$$\{-1\geq\cdots\geq-(\nu-2),-\alpha,-\beta\}=$$
$$\{-1+(1-q)\geq\cdots\geq-(\nu-2)+(1-q),\alpha-(\nu-2)-q,\beta-(\nu-2)-q\},$$
since
$P(t)=(-1)^\nu P(-t-(\nu-2)-q)$.
By adding the both sides of the last equality,
$$\alpha+\beta=(\nu/2)(q+1)-1.$$
Hence,
$\coh{0}{F}{A_F}=(d/2)(q+1)+(\nu-1)$.
Formally, we have 
$$\Delta(F,A_F )=(d/2)(1-q)+1.\eqno(\ref{prfo}.1)$$
From \cite[Lemma~4.3]{na:3adc}, we have
$$\dsh{F}\cong-(\nu-2+q)A_F .\eqno(\ref{prfo}.2)$$
From [ibid., Lemma~5.1], we have three cases,
$q=1,2$
or
$3$.

\medskip

{\bf Case 1) $q=1$:}
From [ibid., Lemmas~4.6 and 4.2],
$Y$
has $1$-factorial, terminal singularities, and 
$\rest{\nbl{E/X}}{F}\cong-A_F$.
If we assume the irreducibility of
$F$, 
then from the same argument as in Case~2) in Section~\ref{prfz}, it follows that
$(F,A_F )$
is a Del Pezzo variety.
On the one hand, F2 in Lemma~\ref{mlm} shows the irreducibility of
$F$
for
$n=\dim X=5$.
Therefore, we obtain the statement in Theorem~\ref{fmth}.

\smallskip

\noindent
Note that F1 in Lemma~\ref{mlm} guarantees that
$F$
is irreducible for Cases 2) and 3) below.

\medskip

{\bf Case 2) $q=2$:}
As above
$Y$
has $2$-factorial, terminal singularities, and 
$\rest{\nbl{E/X}}{F}\cong-2A_F$.
The fact that
$F$
is irreducible enables us to use the non-negativity of the delta genus.
Hence, from (\ref{prfo}.1), we deduce that
$d=2$
and thus
$\Delta(F,A_F )=0$.
Hence, from \cite[Corollary~(4.3)]{fj:sdt}, we obtain that
$(F,A_F ,\rest{\nbl{E/X}}{F})\cong(\hyq{n-2} ,\str(1),\str(-2))$.

\medskip

{\bf Case 3) $q=3$:}
As above
$Y$
has $3$-factorial, terminal singularities, and 
$\rest{\nbl{E/X}}{F}\cong-3A_F$.
From (\ref{prfo}.1),
$d=1$
and thus
$\Delta(F,A_F )=0$.
Hence,
$(F,A_F ,\rest{\nbl{E/X}}{F})\cong(\prs{n-2} ,\str(1),\str(-3))$.
As in the Proof of \cite[Theorem~4]{fj:adj}, we can improve the description on a fiber
$F$.

\medskip

\section{Proof of the fourth adjoint contraction of divisorial type with
$\protect{\dim}f(E)=2$}
\label{prft}

\noindent
Let
$F$
be a general fiber of
$E\to f(E)$,
$\nu=\dim F \;(=n-3\geq2),\;d={A_F}^\nu$
and
$P(t)=\chi(F,tA_F )$.
From \cite[Lemma~4.5]{na:3adc}, as before, we have

$$P(t)=(d/\nu!)(t+1)\cdots(t+\nu-1)(t+\alpha)$$
for some complex number
$\alpha$,
$$\nu=d\alpha,$$
$$\coh{0}{F}{A_F}=d+\nu,$$
and
$$\{-1\geq\cdots\geq-(\nu-1),-\alpha\}=$$
$$\{-1+(1-q)\geq\cdots\geq-(\nu-1)+(1-q),\alpha-(\nu-1)-q\},$$
since
$$P(t)=(-1)^\nu P(-t-(\nu-1)-q).$$
By adding the both sides,
$$\alpha=q\nu/2,$$
and thus
$\nu=qd\nu/2$,
or
$2=qd$.
Hence,
$$(q,d)=(2,1)\mbox{ or }(1,2).$$

\medskip

\noindent
Note that from F1 in Lemma~\ref{mlm},
$F$
is irreducible for Cases 1) and 2) below.

\medskip

{\bf Case 1)}
$(q,d)=(2,1)$:
$Y$
has $2$-factorial, terminal singularities, and
$\rest{\nbl{E/X}}{F}\cong-2A_F$.
Since
$\Delta(F,A_F )=0$,
we have
$(F,A_F ,\rest{\nbl{E/X}}{F})\cong(\prs{n-3},\str(1),\str(-2))$.

\medskip

{\bf Case 2)}
$(q,d)=(1,2)$:
$Y$
has $1$-factorial, terminal singularities, and
$\rest{\nbl{E/X}}{F}\cong-A_F$.
From
$\Delta(F,A_F )=0$,
we deduce
$(F,A_F ,\rest{\nbl{E/X}}{F})\cong(\hyq{n-3},\str(1),\str(-1))$.

\medskip

\section{Proof of the fourth adjoint contraction of divisorial type with
$\protect{\dim}f(E)=3$}
\label{prftd}

\noindent
Let
$F$
be a general fiber of
$E\to f(E)$,
$\nu=\dim F \;(=n-4\geq1),\;d={A_F}^\nu$
and
$P(t)=\chi(F,tA_F )$.
As above, we have
$$P(t)=(d/\nu!)(t+1)\cdots(t+\nu),$$
$d=1,\; \coh{0}{F}{A_F}=\nu+1$,
and
$$\{-1\geq\cdots\geq-\nu\}=\{-1+(1-q)\geq\cdots\geq-\nu+(1-q)\},$$
since
$P(t)=(-1)^\nu P(-t-\nu-q)$.
From these, we deduce that
$F$
is irreducible,
$\Delta(F,A_F )=0$
and
$q=1$.
Hence,
$Y$
has $1$-factorial, terminal singularities, and
$(F,A_F ,\rest{\nbl{E/X}}{F})\cong(\prs{n-4},\str(1),\str(-1))$.

\medskip

\section{Proof of the irreducibleness of a general fiber $F$ of $E\to f(E)$}
\label{pirr}
\noindent
First of all, we will make clear the situation to deal with.
Let
$f:X\to Y$
be the $4$-th adjoint contraction of divisorial type with
$r:=\dim f(E)$,
and let
$E$
be a unique exceptional divisor for 
$f$.
Recall \cite[Lemma~4.3]{na:3adc}:
For general hyperplane sections
$V_1 ,\ldots\hspace{-.001in},V_r$
on
$Y$,
we define
$M_i =f^\ast V_i$, 
a total transform of 
$V_i$,
$$M=M_1 \cap\cdots\cap M_r
\mbox{ and } 
V=V_1 \cap\cdots\cap V_r.$$
Let
$F$
be a general fiber of 
$\rest{f}{E} :E\to f(E)$.
Note that
$F$
can be considered as a connected component of
$M\cap E$.
Let
$$\nu:=\dim F =n-1-r \:(n=\dim X).$$
A positive integer
$q$
is defined by 
$$\rest{E}{F}\cong-q A_F.$$

\noindent
Note also that
$F$
is a connected projective reduced scheme of equidimension
$\nu$
and a local complete intersection in
$X$.
In order to complete the arguments in the previous sections, we have to show

\begin{LE}
\label{mlm}
\begin{description}
\item[F1]
$F$
is irreducible for all 
$(r,q)=(1,2), (1,3), (2,1).$

\item[F2]
$F$
is irreducible for
$(r,q)=(1,1)$ if
$n=5$
(or equivalently,
$\nu=3$).
\end{description}
\end{LE}
\noindent
First of all, we make the following
\begin{AS}
\label{redass}
$F$
is not irreducible.
\end{AS}

\noindent
In what follows, this assumption will be kept until we obtain a contradiction at the end of this section.
Note that the scheme structure of
$F$
requires
$F$
to be an integral divisor on the smooth variety
$M$
and to have the decomposition into distinct prime divisors
$G_i$'s
on
$M$:
$$F=\sum_{i\geq1} G_i .$$
Let
$G$
be an arbitrary irreducible component of
$F$,
say
$G=G_1$,
and
$G'=\sum_{i\geq2} G_i .$
Restrict
$F=G+G'$
to
$G$,
and
$$\rest{F}{G}=\nbl{G/M} +D,$$
where
$D:=\rest{G'}{G}$.
Note also that the Assumption~\ref{redass} forces
$G'$
to be a nonzero effective Cartier divisor on
$M$,
and thus the connectedness of
$F$
implies that
$D$
is a nonzero effective Cartier divisor on
$G$,
or
$D>0$.
On the one hand, \cite[Lemma~4.2~(2)]{na:3adc} implies
$$\rest{F}{G}\cong\rest{(\rest{F}{F})}{G}\cong\rest{\nbl{F/M}}{G}\cong\rest{(\rest{\nbl{E/X}}{F})}{G}\cong-qA_G . \eqno(\ref{pirr}.1)$$
Therefore, altogether we obtain
$$-qA_G -\nbl{G/M}=D>0 \mbox{ in }\Div(G), \: \eqno(\ref{pirr}.2)$$
where
$\Div(G)$
denotes the set of all Cartier divisors on
$G$.
From the adjunction formula,
$$\begin{array}{lcl}
\dsh{G}&\cong&\rest{(K_M +G)}{G} \\
&\cong&\rest{[\rest{(K_X +M)}{M} +G]}{G} \\
&=&\rest{K_X}{G}+\rest{G}{G}.
\end{array}$$
Apply [ibid., Lemma~4.2~(1)] to this, and from (\ref{pirr}.2),
$$\begin{array}{rcl}
\quad\quad\quad\quad\quad\quad\quad\quad\quad\quad\quad\quad\quad
\dsh{G}&=&-(\nu+r-3)A_G +\nbl{G/M}\\
       &=&-(\nu+r+q-3)A_G -D 
\quad\quad\quad\quad\quad\quad\quad\quad\quad\quad\quad\quad\quad(\ref{pirr}.3)
\end{array}$$

\noindent
Let
$g:S\to G$
be a desingularization of 
$G$,
and
$A_S :=g^* A_G$.
Then, from [ibid., (4) in Proof of Lemma~5.1],
$$(g_* K_S )\cdot {A_G}^{\nu-1} \geq K_S \cdot {A_S}^{\nu-1} \eqno(\ref{pirr}.4)$$
The same argument as in [ibid., (1), (2) and (3) in Proof of Lemma~5.1] can apply to
$\dsh{G}$
and results in
$$\dsh{G}\cdot{A_G}^{\nu-1}\geq g_* K_S\cdot{A_G}^{\nu-1} \eqno(\ref{pirr}.5)$$

\noindent
Let
$\mu:\tilde{G}\to G$
be the normalization of an arbitrary irreducible component
$G$
of
$F$,
and
$A_{\tilde{G}} =\mu^* A_G$.

\begin{LE}
\label{clofg}
If
$(r,q)=(1,2),(1,3)$
or
$(2,1)$,
then
$$(\tilde{G},A_{\tilde{G}} )\cong (\prs{\nu}, \str(1)).$$

\noindent
If
$(r,q)=(1,1)$
and if
$n=\dim X=5$
or
$\nu=3$,
then
$$(\tilde{G},A_{\tilde{G}} )\cong(\prs{\nu}, \str(1))\mbox{ or }(\hyq{\nu},\str(1)).$$
\end{LE}

\begin{pf}
The proof is rather lengthy.
We assume that 
$(\tilde{G},A_{\tilde{G}} )\not\cong (\prs{\nu}, \str(1))$.
Then, \cite[Theorem(2.2)]{fj:rqp} implies one of the following two conditions.

$$\mbox{There exists an integer
$t$
with
$1\leq t<\nu$
such that
$\Coh{0}{S}{K_S +tA_S}\not=0$.} \eqno(\ref{pirr}.a)$$
$$\mbox{$\Coh{0}{S}{K_S +\nu A_S}\not=0$
and
$\Coh{0}{S}{K_S +tA_S}=0$
for
any
$1\leq t<\nu$.}  \eqno(\ref{pirr}.b)$$

\begin{CL}
\label{nota}
Case
$(\ref{pirr}.a)$
does not occur.
\end{CL}

\begin{pf}
If Case 
$(\ref{pirr}.a)$ 
did occur,
$(K_S +(\nu-1) A_S )\cdot{A_S}^{\nu-1} \geq0$.
Hence, from (\ref{pirr}.4), (\ref{pirr}.5) and (\ref{pirr}.3), it follows that

$$\begin{array}{rcl}
(\nu-1) {A_S}^\nu
&\geq&-K_S \cdot{A_S}^{\nu-1}\\
&\geq&-(g_\ast K_S )\cdot {A_G}^{\nu-1}\\
&\geq&-\dsh{G}\cdot {A_G}^{\nu-1}\\
&=&(\nu+r+q-3) {A_G}^\nu +D\cdot {A_G}^{\nu-1} \\
&>&(\nu+r+q-3) {A_G}^\nu ,
\end{array}$$
where we used
$D>0$
 for the last inequality.
Thus,
$r+q<2$,
a contradiction.
\end{pf}

\noindent
Therefore, we are in Case 
$(\ref{pirr}.b).$
Since 
$\Coh{0}{S}{K_S +\nu A_S}\not=0$,
we have
$(K_S +\nu A_S )\cdot{A_S}^{\nu-1} \geq0$.
Hence, exactly the same argument as in the Proof of Claim~\ref{nota} implies that
$$\nu {A_S}^\nu >(\nu+r+q-3) {A_G}^\nu$$
and thus
$$r+q<3\eqno(\ref{pirr}.6).$$

\smallskip

\noindent
If
$(r,q)=(1,2), (1,3)$
or
$(2,1)$, 
then (\ref{pirr}.6) implies a contradiction.
Therefore, 
$(\tilde{G},A_{\tilde{G}} )\cong (\prs{\nu}, \str(1))$
for these three cases, and the first statement of Lemma~\ref{clofg} is proven.

\smallskip

\noindent
From now on, we will prove the second statement of Lemma~\ref{clofg}.
To prove this, it suffices to show that
$(\tilde{G},A_{\tilde{G}} )\cong(\hyq{\nu},\str(1))$
under the conditions
$(r,q)=(1,1)$
and 
$\nu=3$
as well as
$(\ref{pirr}.b)$.

\begin{LE}
\label{fstsdg}
If
$\coh{0}{S}{K_S +\nu A_S }\not=1$,
then
$$g(G,A_G )=g(\tilde{G},A_{\tilde{G}} )=\Delta(\tilde{G},A_{\tilde{G}} )=0,$$
and
${A_{\tilde{G}}}^\nu >2$.
\end{LE}

\begin{pf}
From \cite[Corollary~(2.8)]{fj:rqp},
$\coh{0}{S}{K_S +\nu A_S }\not=1$
implies
$$(K_S +\nu A_S )\cdot{A_S}^{\nu-1} >0.$$
As in the Proof of Claim~\ref{nota},
$$\begin{array}{rcl}
\nu {A_S}^\nu
&>&-K_S \cdot{A_S}^{\nu-1}\\
&\geq&-\dsh{G}\cdot {A_G}^{\nu-1}\\
&>&(\nu-1) {A_G}^\nu.
\end{array}$$
Thus,
$$-{A_G}^\nu <(\dsh{G}+(\nu-1)A_G )\cdot{A_G}^{\nu-1} <0.$$
Note that the scheme structure of
$F$
enforces
$(G,A_G )$
to be a locally Gorenstein pre-polarized variety (Definition~\ref{pv}).
Hence, we are able to apply Proposition~\ref{sgf} in order to have
$-{A_G}^\nu <2g(G,A_G )-2<0,$
or
$$1-\frac{{A_G}^\nu}{2}<g(G,A_G )<1. \eqno(\ref{pirr}.7)$$
In particular,
${A_{\tilde{G}}}^\nu ={A_G}^\nu >2$.
On the one hand, from Proposition~\ref{inse}, it follows that
$$g(G,A_G )\geq g(\tilde{G},A_{\tilde{G}} ).$$
On the other hand, the restriction that
$\nu=3$
enables us to apply Lemma~\ref{pssg} (\cite[Corollary~(4.8)]{fj:rqp}) to a normal polarized variety
$(\tilde{G},A_{\tilde{G}} )$
of dimension three which asserts that
$$g(\tilde{G},A_{\tilde{G}} )\geq0,$$
and that if
$g(\tilde{G},A_{\tilde{G}} )=0$,
then
$\Delta(\tilde{G},A_{\tilde{G}} )=0$.
Therefore, (\ref{pirr}.7) will yield the desired result in Lemma~\ref{fstsdg}.
\end{pf}

\noindent
Let
$\tilde{D}=\mu^* D$,
the pull-back of the positive Cartier divisor
$D$.
Then,
$\tilde{D}$
is also a positive Cartier divisor on
$\tilde{G}$
since locally,
$\mu_* \str_{\tilde{G}}$
is the integral closure of an integral domain
$\str_G$.

\begin{LE}
\label{gnml}
If
$g(G,A_G )=g(\tilde{G},A_{\tilde{G}} )=0$,
then
$$\begin{array}{rcl}
K_{\tilde{G}}
&\cong&\mu^* \dsh{G} \\
&\cong&-(\nu-1)A_{\tilde{G}} -\tilde{D},
\end{array}$$
and
$$\tilde{D}\cdot{A_{\tilde{G}}}^{\nu-1} =2.$$
\end{LE}

\begin{pf}
From Proposition~\ref{inse},
$g(G,A_G )=g(\tilde{G},A_{\tilde{G}} )$
implies that the non-normal locus
$G_{\mbox{non}}$
of
$G$
has dimension less than
$\nu-1$,
or
$\dim G_{\mbox{non}} <\nu-1$.
On the one hand,
$$\rest{\mu}{\tilde{G}\setminus\mu^{-1} G_{\mbox{non}}}:\tilde{G}\setminus\mu^{-1} G_{\mbox{non}} \to G\setminus G_{\mbox{non}}$$
is an isomorphism.
Hence,
$K_{\tilde{G}} \cong\mu^* \dsh{G}$
on
$\tilde{G}\setminus\mu^{-1} G_{\mbox{non}}$.
However, this isomorphism extends to the whole normal variety
$\tilde{G}$
since
$\dim \mu^{-1} G_{\mbox{non}} <\nu-1$,
and
$\mu^* \dsh{G}$
is a reflexive sheaf of rank one on
$\tilde{G}$.
The rest of the first assertion follows from (\ref{pirr}.3) and
$(r,q)=(1,1)$.

\smallskip

\noindent
As for the second assertion, combine the assumption with \cite[Lemma~(1.8)]{fj:rqp}, and 
$$g(S,A_S )=g(\tilde{G},A_{\tilde{G}} )=0.$$
Hence, we obtain
$$(K_S +(\nu-1) A_S )\cdot{A_S}^{\nu-1} =2g(S,A_S )-2=-2.$$
This equality, (\ref{pirr}.4), (\ref{pirr}.5) and (\ref{pirr}.3) imply

$$\begin{array}{rcl}
(\nu-1) {A_S}^\nu
&=&-K_S \cdot{A_S}^{\nu-1} -2\\
&\geq&-(g_\ast K_S )\cdot {A_G}^{\nu-1} -2\\
&\geq&-\dsh{G}\cdot {A_G}^{\nu-1} -2\\
&=&(\nu-1) {A_G}^\nu +D\cdot {A_G}^{\nu-1} -2\\
&\geq&(\nu-1) {A_G}^\nu -1.
\end{array}$$
Therefore,
$$-\dsh{G}\cdot {A_G}^{\nu-1} -2=(\nu-1) {A_G}^\nu -1 
\mbox{ or } (\nu-1) {A_G}^\nu.$$
Note that the first possibility is easily ruled out.
Otherwise, from Proposition~\ref{sgf},
$$\begin{array}{rcl}
2g(G,A_G )-2
&=&(\dsh{G}+(\nu-1)A_G )\cdot{A_G}^{\nu-1} \\
&=&-1.
\end{array}$$
This is impossible since
$g(G,A_G )\in\integer$. \\
 \\
From the second case, it follows that
$$D\cdot{A_G}^{\nu-1} =2,$$
and thus
$\tilde{D}\cdot{A_{\tilde{G}}}^{\nu-1} =2$.
\end{pf}

\begin{LE}
\label{fstask}
Assume that the condition 
$(\ref{pirr}.b)$
holds,
$(r,q)=(1,1)$
and 
$\nu=3$.
Then,
$\coh{0}{S}{K_S +\nu A_S }=1$.
\end{LE}

\begin{pf}
If
$\coh{0}{S}{K_S +\nu A_S }\not=1$,
then from the previous Lemmas~\ref{fstsdg} and \ref{gnml}, it would follow that
$$\begin{array}{rcl}
g(\tilde{G},A_{\tilde{G}} )
&=&\Delta(\tilde{G},A_{\tilde{G}} )=0, \\
{A_{\tilde{G}}}^\nu
&>&\tilde{D}\cdot {A_{\tilde{G}}}^{\nu-1} =2, \mbox{ and }\\
K_{\tilde{G}}
&\cong&-(\nu-1)A_{\tilde{G}} -\tilde{D} \quad (\tilde{D}>0,\: \nu=3).
\end{array}$$
From
$\Delta(\tilde{G},A_{\tilde{G}} )=0$,
\cite[Theorem~(4.11)]{fj:sdt} will provide the structure of the normal polarized variety
$(\tilde{G},A_{\tilde{G}} )$
as follows.
Since
$A_{\tilde{G}}$
is very ample,
$\lsys{A_{\tilde{G}}}$
induces the embedding
$$\tilde{G}\hookrightarrow \prs{N},$$
where
$N={A_{\tilde{G}}}^\nu +\nu-1$,
and we identify 
$\tilde{G}$
with its image in 
$\prs{N}$.
Note that 
$$A_{\tilde{G}} \cong\rest{\str_{\prs{N}} (1)}{\tilde{G}}, \mbox{ and } \dim \tilde{G}_\sing \leq\nu-2,$$
where
$\tilde{G}_\sing$
is the singular locus of
$\tilde{G}$.
From [ibid.], there is a linear subspace
$T$
of
$\prs{N}$
such that
\begin{enumerate}
\item
$\codim(T,\prs{N})=\dim \tilde{G}_\sing +1$,
\item
$M:=\tilde{G}\cap T$
is a smooth subvariety of 
$\tilde{G}$
with
$\Delta(M,A_M )=0$,
where
$A_M :=\rest{A_{\tilde{G}}}{M}$,
and
\item
$\tilde{G}=M\ast \tilde{G}_\sing$,
the union of lines each of which goes through two different points, one in
$M$
and the other in
$\tilde{G}_\sing$.
\end{enumerate}

\noindent
We note that if 
$m=\dim M<\dim \tilde{G}=\nu$,
then
$\tilde{G}$
must be singular.
Therefore, from \cite[Corollary~(4.13)]{fj:sdt}, it follows that
$$\Pic(\tilde{G})\cong\integer\cdot A_{\tilde{G}} .$$
Hence, there is a positive integer
$l$
such that
$\tilde{D}\sim lA_{\tilde{G}}$,
and thus
$$\tilde{D}_M :=\rest{\tilde{D}}{M}\sim lA_M \quad(0<l\in\integer). \eqno(\ref{pirr}.8)$$
From [ibid., Theorems~(4.9) and (4.11)], the condition that
${A_{\tilde{G}}}^\nu \geq3$
possibly allows
$(M,A_M )$
to be only either a Veronese surface or a scroll over a $1$-dimensional projective space.
We will show that neither of these two cases occurs, which will, in turn, show
$\coh{0}{S}{K_S +\nu A_S }=1$
and finish the proof of Lemma~\ref{fstask}.

\begin{CL}
\label{vrnon}
$(M,A_M )$
cannot be a Veronese surface
$(\prs{2},\str(2))$.
\end{CL}

\begin{pf}
If this were in case, it would follow that
$(M,A_M )\cong(\prs{2},\str(2))$,
and thus
$m=2<3=\nu$.
Hence, from (\ref{pirr}.8),
$\tilde{D}_M \cong\str_{\prs{2}} (2l)$.
Therefore,
$$\begin{array}{rcl}
\tilde{D}\cdot {A_{\tilde{G}}}^{\nu-1}
&=&\tilde{D}_M \cdot {A_M}^{\nu-2}\\
&=&2l\cdot2\\
&=&4l.
\end{array}$$
This contradicts that
$\tilde{D}\cdot {A_{\tilde{G}}}^{\nu-1} =2$.
\renewcommand{\qed}{\hfill{q.e.d.~of~Claim~\ref{vrnon}.}}
\end{pf}

\begin{CL}
\label{scnon}
$(M,A_M )$
cannot be a scroll.
\end{CL}

\begin{pf}
If so, it would follow that
$$(M,A_M )\cong({\Bbb P}({\cal E}),\xi)$$
which has a
$\prs{m-1}$-bundle structure
$$\pi:M\longrightarrow \prs{1},$$
where
$m=\dim M$,
$${\cal E}\cong\str_{\prs{1}}(a_1 )\oplus\cdots\oplus\str_{\prs{1}}(a_m ),$$
and
$a_i >0$
for any $i$.
Since 
$\codim(T,\prs{N})\leq\nu-1$,
we have
$m=\dim M\geq\nu-(\nu-1)=1$,
and thus
$$m=1,2
\mbox{ or }
3=\nu.$$

\begin{description}
\item[Case $m=1$]
$\coh{1}{M}{\str_M}=g(M,A_M )=0$
since
$g(\tilde{G},A_{\tilde{G}} )=0$.
Hence,
$M\cong\prs{1}$.
On the one hand, 
$$\deg A_M ={A_{\tilde{G}}}^\nu \geq3,$$
and thus
$A_M \cong\str_{\prs{1}} (d), \; (d\geq3)$.
Since
$m=1<3=\nu$,
from (\ref{pirr}.8),
$\tilde{D}_M \cong\str_{\prs{1}} (ld)$.
Hence,
$$\tilde{D}\cdot {A_{\tilde{G}}}^{\nu-1} =\deg \tilde{D}_M =ld\geq3,$$
a contradiction to
$\tilde{D}\cdot {A_{\tilde{G}}}^{\nu-1} =2$.

\smallskip

\item[Case $m\geq2$]
We do not use (\ref{pirr}.8), since this case includes the one that
$\tilde{G}=M$
or
$\tilde{G}$
is nonsingular, and \cite[Corollary~(4.13)]{fj:sdt} is not applicable.
From \cite[II, Exercise~7.9]{ha:alg},
$$\Pic(M)\cong\Pic(\prs{1})\times\integer\cong\integer{\bf f}+\integer A_M,$$
where
${\bf f}$
is a fiber of 
$\pi$.
From the adjunction formula,
$K_M \cong\rest{[K_{\tilde{G}} +(\nu-m)A_{\tilde{G}} ]}{M}$.
Hence,
$$\begin{array}{rcl}
K_M &\cong&\rest{[-(\nu-1)A_{\tilde{G}} -\tilde{D}+(\nu-m)A_{\tilde{G}} ]}{M} \\
  &=&-(m-1)A_M -\tilde{D}_M.
\end{array}$$
Since
$\tilde{D}>0$,
$\tilde{D}_M \cong a{\bf f}+bA_M$
for some integers
$a,b\geq0$.
Hence,
$$K_M \cong-a{\bf f}-(m+b-1)A_M. \eqno(\ref{pirr}.9)$$
From the canonical bundle formula,
$$\begin{array}{rcl}
\quad\quad\quad\quad\quad\quad\quad\quad\quad\quad\quad\quad
K_M &\cong&\pi^* (K_{\prs{1}} +\det{\cal E})-m\xi \\
  &\cong&\pi^* (\str_{\prs{1}} (-2)+\str_{\prs{1}} (\sum_{i=1}^m a_i ))-mA_M \quad\quad\quad\quad\quad\quad\quad\quad\quad\quad\quad\\
  &\cong&(\sum_{i=1}^m a_i -2){\bf f}-mA_M. \hfill(\ref{pirr}.10)
\end{array}$$
Compare (\ref{pirr}.10) with (\ref{pirr}.9), and
$$\sum_{i=1}^m a_i +a=2.$$
This is absurd since
$$\sum_{i=1}^m a_i =\deg(\det{\cal E})={A_M}^m \geq3.$$
Hence, Case
$m\geq2$
does not occur.
\end{description}
\renewcommand{\qed}{\hfill{q.e.d.~of~Claim~\ref{scnon}.}}
\end{pf}

\noindent
From Claims~\ref{vrnon} and \ref{scnon}, we have proven that 
$\coh{0}{S}{K_S +\nu A_S }=1$.
\renewcommand{\qed}{\hfill{q.e.d.~of~Lemma~\ref{fstask}.}}
\end{pf}

\medskip

\noindent
Recall that we are in Case 
$(\ref{pirr}.b)$
in the Proof of the second statement of Lemma~\ref{clofg}.
From this and Lemma~\ref{fstask},
the conditions which we have obtained are the following.
$$\Coh{0}{S}{K_S +tA_S}=0
\mbox{ for any }
1\leq t<\nu,
\mbox{ and }
\coh{0}{S}{K_S +\nu A_S }=1.$$
Thus, \cite[Theorem~(2.3)]{fj:rqp} implies either
\begin{enumerate}
\item
${A_{\tilde{G}}}^\nu =g(\tilde{G},A_{\tilde{G}} )=1$,
or
\item
$(\tilde{G},A_{\tilde{G}} )\cong(\hyq{\nu},\str(1))$.
\end{enumerate}

\noindent
If the first case occurred, it would follow from [ibid., Lemma~(1.8)] that
$$\begin{array}{rcl}
(K_S +(\nu-1)A_S )\cdot{A_S}^{\nu-1} &=& 2g(S,A_S )-2\\
                                     &=& 2g(\tilde{G},A_{\tilde{G}} )-2 \\
                                     &=& 0
\end{array}$$
Hence,
$(\nu-1){A_S}^\nu =-K_S \cdot{A_S}^{\nu-1}$,
which leads to a contradiction exactly as in the Proof of Claim~\ref{nota}.
Therefore, only the second case above will possibly occur, and we have proven Lemma~\ref{clofg}.
\end{pf}

\medskip

\noindent
We can express
$\nbl{G/M}$
and
$D$
by
$A_G$
so as to get a more detailed list for
$(\tilde{G},A_{\tilde{G}} )$.

\begin{LE}
\label{clofgand}
The normalization
$\mu:\tilde{G}\to G$
of an arbitrary irreducible component
$G$
of
$F$
satisfies the following list.

\begin{equation*}
(\tilde{G},A_{\tilde{G}} ,\nbl{G/M},D)\cong
\begin{cases}
(\prs{\nu},\str_{\prs{\nu}} (1),-3A_G ,0A_G )&\text{if $(r,q)=(1,3)$} \hspace{4cm}\\
(\prs{\nu},\str_{\prs{\nu}} (1),-(4-r)A_G ,A_G )&\text{if $(r,q)=(1,2)$ or $(2,1)$} \\
(\prs{\nu},\str_{\prs{\nu}} (1),-3A_G ,2A_G ),&\text{ } \\
\mbox{or}                       &\text{if $(r,q)=(1,1)$ and if $\nu=3$} \\
(\hyq{\nu},\str_{\hyq{\nu}} (1),-2A_G ,A_G )&\text{ }
\end{cases}
\end{equation*}
\end{LE}

\begin{pf}
Either case in Lemma~\ref{clofg} implies that
$g(\tilde{G},A_{\tilde{G}} )=0$.
On the other hand, from Proposition~\ref{sgf} and (\ref{pirr}.3),
$$\begin{array}{rcl}
2g(G,A_G )-2
&=&(\dsh{G}+(\nu-1)A_G )\cdot{A_G}^{\nu-1} \\
&=&-(r+q-2){A_G}^{\nu} -D\cdot{A_G}^{\nu-1} \\
&<&0,
\end{array}$$
and thus
$g(G,A_G )\leq0$.
From Proposition~\ref{inse},
$g(G,A_G )\geq g(\tilde{G},A_{\tilde{G}} ).$
Hence,
$g(G,A_G )=g(\tilde{G},A_{\tilde{G}} )=0$,
and from Proposition~\ref{inse}, this implies
$$\dim G_{\mbox{non}} <\nu-1. \eqno(\ref{pirr}.c)$$
Therefore, the same argument as in Lemma~\ref{gnml} shows that
$$\begin{array}{rcl}
K_{\tilde{G}}
&\cong&\mu^* \dsh{G} \\
&\cong&-(\nu+r+q-3)A_{\tilde{G}} -\tilde{D}.
\end{array}$$
Hence, from Lemma~\ref{clofg}, it follows that
\begin{equation*}
\tilde{D}\cong
\begin{cases}
0A_{\tilde{G}}  &\text{if $(r,q)=(1,3)$} \hspace{4cm}\\
A_{\tilde{G}}   &\text{if $(r,q)=(1,2)$ or $(2,1)$} \\
2A_{\tilde{G}} ,&\text{ } \\
\mbox{or}       &\text{if $(r,q)=(1,1)$ and if $\nu=3$} \\
A_{\tilde{G}}   &\text{ }
\end{cases}
\end{equation*}
To obtain the result on
$D$
in the statement of Lemma~\ref{clofgand}, it suffices to show the following

\begin{CL}
\label{divrl}
Let
$B$
be any Cartier divisor on
$G$.
If
$\mu^* B\cong\str_{\tilde{G}}$,
then
$B\cong\str_G$.
\end{CL}

\begin{pf}
We apply \cite[p.65, Proposition]{mu:cas}:
\begin{FA}
\label{mmecc}
For two Cartier divisors
$D_1$
and
$D_2$
on a Noetherian scheme
$X$,
$D_1 =D_2$
as divisors if and only if the images of
$D_1$
and
$D_2$
in the stalk
$({\cal K}_X ^* /\str_X ^* )_p$
are equal for all
$p\in X$
whenever
$\depth(\str_{X,p} )=1$,
where
${\cal K}_X$
denotes the sheaf of total quotient rings of
$X$.
\end{FA}

\noindent
Since
$\mu^* B\cong\str_{\tilde{G}}$,
there is a rational function
$\tilde{f}\in\Coh{0}{\tilde{G}}{{\cal K}_{\tilde{G}} ^*}$
such that
$\mu^* B=\pdiv(\tilde{f})$,
where
$\pdiv:\Coh{0}{\tilde{G}}{{\cal K}_{\tilde{G}} ^*}\to
\Coh{0}{\tilde{G}}{{\cal K}_{\tilde{G}} ^* /\str_{\tilde{G}} ^*}$
is the natural homomorphism.
Since
$\mu:\tilde{G}\to G$
is birational,
$\mu^* :\Coh{0}{G}{{\cal K}_G ^*}\to\Coh{0}{\tilde{G}}{{\cal K}_{\tilde{G}} ^*}$
is an isomorphism.
Hence, there is a rational function
$f\in\Coh{0}{G}{{\cal K}_G ^*}$
such that
$\mu^* (f)=\tilde{f}$,
and thus
$\mu^* B=\mu^* \pdiv(f)$.
Note that through the isomorphism
$\rest{
\mu
}{
\tilde{G}\setminus\mu^{-1} G_{\non}
}
:\tilde{G}\setminus\mu^{-1} G_{\mbox{non}} \to G\setminus G_{\mbox{non}}$,
it follows that
$B=\pdiv(f)$
on
$G\setminus G_{\mbox{non}}$.
In other words, the condition in the above Fact~\ref{mmecc} is satisfied except on the non-normal locus
$G_{\mbox{non}}$
of
$G$.
However, we can show that
$$\depth\str_{G,p} \geq2\mbox{ for any } p\in G_{\non}.$$
Because
$G$
is locally Cohen-Macaulay,
$\str_{G,p}$
is a Cohen-Macaulay local ring, and hence
$$\begin{array}{rcl}
\depth \str_{G,p}
&=&\dim\str_{G,p}\\
&\geq&2,
\end{array}$$
where we applied (\ref{pirr}.c) to the last inequality.
Therefore, from Fact~\ref{mmecc}, we derive that
$B=\pdiv(f)$
on the whole
$G$,
or
$B\cong\str_G$.
\renewcommand{\qed}{\hfill{q.e.d.~of~Claim~\ref{divrl}.}}
\end{pf}
We simply apply (\ref{pirr}.2) to calculate
$\nbl{G/M}$
based on the above results on
$D$.
\renewcommand{\qed}{\hfill{q.e.d.~of~Lemma~\ref{clofgand}.}}
\end{pf}

\noindent
Now we can finish up the Proof of Lemma~\ref{mlm}.
We will examine the results of Lemma~\ref{clofgand} in three cases according that
$q=3,2$
or
$1$.
Recall that the irreducible components
$G$'s
of
$F$
are prime Cartier divisors on the smooth variety
$M$,
and each of the components intersects at least one of the others.
\begin{description}
\item[Case $q=3$]
From (\ref{pirr}.2),
$D>0$,
a contradiction.

\item[Case $q=2$]
$\tilde{D}\cong A_{\tilde{G}} \cong\str_{\prs{\nu}} (1)$.
This implies that
$D$
is a prime divisor.
Hence, each component intersects with only one of the other components.
Since
$F$
is connected, there are only two irreducible components total:
$F=G+G'$.
Thus, taking an integral curve
$Z$
in
$G\cap G'$,
we have
$$\begin{array}{rcl}
F\cdot Z&=&G\cdot Z+G'\cdot Z \\
        &=&\nbl{G/M}\cdot Z+\nbl{G'/M}\cdot Z \\
        &=&-3A_G \cdot Z-3A_{G'} \cdot Z.
\end{array}$$
On the other hand, from (\ref{pirr}.1),
$F\cdot Z=-qA_G \cdot Z=-2A_G \cdot Z.$
This is absurd.

\item[Case $q=1$]
There are possibly three kinds of irreducible components of
$F$:
\begin{equation*}
\quad\quad\quad\quad\quad\quad\quad\quad
(\tilde{G},A_{\tilde{G}} ,\nbl{G/M},D)\cong
\begin{cases}
a)\quad(\prs{\nu},\str_{\prs{\nu}} (1),-2A_G ,A_G )&\text{}\\
b)\quad(\prs{\nu},\str_{\prs{\nu}} (1),-3A_G ,2A_G )&\text{ $(\ref{pirr}.11)$} \\
c)\quad(\hyq{\nu},\str_{\hyq{\nu}} (1),-2A_G ,A_G )\quad\quad\quad\quad\quad\quad\quad\quad&\text{}
\end{cases}
\end{equation*}

\noindent
First of all, we will show the following claim.

\begin{CL}
\label{upth}
$F$
has at most three irreducible components.
\end{CL}

\begin{pf}
Note first that from (\ref{pirr}.11), the information divisor
$D$
allows each irreducible component of
$F$
to intersect at most two of the other irreducible components.
Assume that
$G$,
$G'$
and
$G''$
are distinct irreducible components of
$F$
and that
$G$
intersects with both
$G'$
and
$G''$.
Since
$\tilde{G}$
is either
$\prs{\nu}$
or
$\hyq{\nu}$,
the pull-back of the positive Cartier divisors
$\rest{G'}{G}$
and
$\rest{G''}{G}$
of
$G$
are ample.
Hence,
$\rest{G'}{G}$
and
$\rest{G''}{G}$
are ample since
$\mu$
is a finite morphism.
This means that if there are three distinct irreducible components
$G$,
$G'$
and
$G''$
whose union
$G\cup G'\cup G''$
is connected, then, we have
$$\dim G\cap G'\cap G''\geq1 \eqno(\ref{pirr}.12)$$
If there were the fourth irreducible component
$G'''$
which meets
$G\cup G'\cup G''$,
say
$G'''\cap G\not=\emptyset$,
then, from (\ref{pirr}.12),
$G'''$
would have to meet
$G\cap G'\cap G''$
since
$\rest{G'''}{G}$
is an ample Cartier divisor on
$G$
by the same argument as above.
This contradicts the restriction posed by
$D$
at the beginning of the proof.
\renewcommand{\qed}{\hfill{q.e.d.~of~Claim~\ref{upth}.}}
\end{pf}

\noindent
From the Claim~\ref{upth},
$F=G+G'\mbox{ or }G+G'+G''$.
From (\ref{pirr}.12), we can choose an integral curve
$Z$
in
$G\cap G'$
or
$G\cap G'\cap G''$
for each case, and thus
$$F\cdot Z=G\cdot Z+G'\cdot Z\mbox{ or }G\cdot Z+G'\cdot Z+G''\cdot Z. \eqno(\ref{pirr}.13)$$
Let
$G^{\#}$
stand for
$G$,
$G'$
or
$G''$.
From (\ref{pirr}.11), it follows that
$G^{\#} \cdot Z=-iA_{G^{\#}} \cdot Z$
for
$i=2$
or
$3$.
On the other hand,
from (\ref{pirr}.1),
$F\cdot Z=-A_{G^{\#}} \cdot Z$.
This value
$-A_{G^{\#}} \cdot Z$
is always strictly larger than either value of the right hand side in (\ref{pirr}.13) for any
$i=1,2$.
This is a contradiction, and thus we have excluded all the possible cases arising from the Assumption~\ref{redass} that
$F$
is not irreducible.
\end{description}
Therefore,
$F$
must be irreducible, and the Proof of Lemma~\ref{mlm} is completed.

\medskip

\section{Appendix: Definition of the $k$-th adjoint contraction}
\label{def}

\noindent
In this appendix, we will prove Proposition~\ref{1stcl}, 
which supports the definition~\ref{defadj} of the 
$k$-th adjoint contraction.
Recall, first of all, the definition of numerical reducedness.

\begin{DE}(cf. \cite[Definition~1.1]{na:3adc})
\label{nrds}
Let 
$f:X\to Y$
be the ray contraction for an extremal ray
$R$.
An $f$-ample Cartier divisor
$A$
is called to be {\em numerically reduced} if there is no ample Cartier divisor
$A'$
such that
$$A\cdot C=p A'\cdot C$$
for an integer $p>1$ and a nonzero effective curve $C$ in $R$,
where a curve means an element in $Z_1 (X)$, the free Abelian group generated by one-dimensional closed integral subvarieties of
$X$.
\end{DE}

\noindent
Given a nonzero effective curve $C$ in $R$, we define
$$ M(C)\stackrel{\mbox{def}}{=}\{A'\cdot C|A'\mbox{ is an $f$-ample Cartier divisor on $X$}\}.$$
$A'\cdot C$ are all positive integers and thus there is an $f$-ample Cartier divisor
$A$
on
$X$
such that
$$A\cdot C=\min M(C)\eqno(\ref{def}.1)$$

\begin{LE}(cf. \cite[Lemma~4.1]{na:3adc})
\label{pnrd}
Let
$A$
be a Cartier divisor defined by (\ref{def}.1).
\begin{enumerate}

\item
$A$
is numerically reduced, and thus a numerically reduced divisor always exists for a given ray
$R$.

\item
For any Cartier divisor
$D$
on 
$X$,
there is a unique integer
$\delta$
such that
$$D\equiv \delta A \mod f^* \Div(Y).$$
In particular, a numerically reduced divisor
$A$
is unique modulo
$f^* \Div(Y)$,
and does not depend on the choice of a curve 
$C$
in 
$R$.

\item
We can choose a numerically reduced divisor
$A$
to be ample by adding 
$f^* (m{\cal A})$,
where
${\cal A}$
is an ample Cartier divisor on
$Y$
and 
$m$ 
is a large enough integer
$m>>1$.
\end{enumerate}
\end{LE}

\begin{pf}
1.
If $A$ were not numerically reduced, then there would be an $f$-ample Cartier divisor $A'$, a nonzero effective curve $C'$ in $R$ and an integer $p>1$ such that 
$$A\cdot C'=p A'\cdot C'.$$
Note that there is a positive real number 
$a$ 
such that
$C'\approx aC$
in 
$R$,
where 
$C$ 
is the fixed curve as in 
$M(C)$, 
and 
$\approx$ 
denotes the numerical equivalence for curves on
$X$.
Hence, 
$A\cdot C=p A'\cdot C$,
a contradiction to 
$A\cdot C=\min M(C)$.
Therefore, $A$ must be numerically reduced.

2.
If $D\cdot C<0$, 
we would consider 
$-D$ 
instead of 
$D$.
So, we can assume that 
$D\cdot C\geq0$.
Divide 
$D\cdot C$ 
by 
$A\cdot C$, 
and we obtain unique integers
$\delta$
and 
$r$
such that
$$D\cdot C=\delta A\cdot C+r, \quad 0\leq r <A\cdot C \quad \mbox{ and }\delta\geq0.$$
Then,
$r=0$.
Otherwise,
$(D-\delta A)\cdot C=r>0$,
and thus 
$D-\delta A$
is $f$-ample.
However, 
$(D-\delta A)\cdot C<A\cdot C$,
which contradicts the choice of
$A$
in (\ref{def}.1).
Hence, 
$D\cdot C=\delta A\cdot C$.
From \cite[Lemma~3-2-5]{kmm}, 
$D\equiv \delta A \; \mod f^* \Div(Y)$.

3.
This is nothing but one of characterizations for $A$ to be $f$-ample.
\end{pf}

\begin{LE}
\label{exdef}
Let 
$f:X\to Y$
be the ray contraction for an extremal ray
$R$.
There is a numerically reduced, ample Cartier divisor 
$A$
on 
$X$
which is unique modulo
$f^* \Div(Y)$,
and there is a unique integer
$k$
with
$-1\leq k\leq n-1$
such that
$$R=(K_X +(n-k)A)^\bot \cap \overline{NE}(X). \eqno(\ref{def}.2)$$
\end{LE}

\begin{pf}
Note that from the first statement in Lemma~\ref{pnrd}, we have a numerically reduced divisor $A$ for $f:X\to Y$.
From \cite[Theorem~3-2-1]{kmm}, there is a nef Cartier divisor 
$H$ 
such that
$R=H^\bot \cap \overline{NE}(X)$,
and 
$K_X \cdot C<0$
for any nonzero curve $C$ in $R$.
Hence, 
$$(H-K_X )\cdot C=-K_X \cdot C>0,$$
and thus
$H-K_X$
is $f$-ample.
Therefore, from  the second statement in Lemma~\ref{pnrd}, there exists a unique integer $l\geq1$ such that
$H-K_X \equiv lA \; \mod f^* \Div(Y)$,
or
$$H \equiv K_X +lA \quad \mod f^* \Div(Y). \eqno(\ref{def}.3)$$
To finish the proof, it suffices to show that
$$R=(K_X +l A)^\bot \cap \overline{NE}(X), \eqno(\ref{def}.4)$$
since from \cite[Theorem~4-2-1]{kmm}, we have $l\leq n+1$, and then the number
$k$ will be defined to be 
$n-l$.

\noindent
To show $\subset$ in (\ref{def}.4),
take
$\gamma\in R$.
Since $\gamma$ is the zero-map on $f^* \Div(Y)$,
(\ref{def}.3) implies that 
$$0=H\cdot\gamma=(K_X +l A)\cdot\gamma.$$
Thus, $\gamma\in(K_X +l A)^\bot \cap \overline{NE}(X)$.

\noindent
On the one hand, from (\ref{def}.3), it follows that
$$H +f^* {\cal A}\sim K_X +lA \mbox{ for some } {\cal A}\in\Div(Y),$$
where $\sim$ denotes the linear equivalence.
Take an ample Cartier divisor 
${\cal L}$ on 
$Y$, 
and add 
$f^* (lm{\cal L})$, 
then
$$H +f^* ({\cal A}+lm{\cal L})\sim K_X +l(A+f^* (m{\cal L})).$$
We can choose an integer 
$m$ 
large enough that 
${\cal A}+lm{\cal L}$ 
is ample.
Hence, writing
${\cal A}$
and 
$A$
for 
${\cal A}+lm{\cal L}$
and
$A+f^* (m{\cal L})$,
respectively,
we obtain
$$H +f^* {\cal A}\sim K_X +lA, \eqno(\ref{def}.5)$$
where 
$A$ 
is a numerically reduced ample divisor, and 
$f^* {\cal A}$ 
and 
$K_X +lA$ are nef.
We will apply this to show 
$\supset$ 
in (\ref{def}.4).
Take any non-zero
$\gamma\in(K_X +l A)^\bot \cap \overline{NE}(X)$.
Then, from the Cone Theorem (\cite[Theorem~4-2-1]{kmm}), 
$\gamma\in\overline{NE}(X)$
implies
$$\gamma\approx\gamma^+ +\sum_{j=1}^r a_j C_j,$$
where
$\gamma^+$
is an element in
$\overline{NE}_{K_X} (X)$,
$C_j$
are integral curves which generate different extremal rays,
$a_j$ 
are non-negative real numbers.
Since 
$\gamma\in(K_X +l A)^\bot$,
$$0=(K_X +l A)\cdot\gamma=(K_X +l A)\cdot\gamma^+ +\sum_{j=1}^r a_j (K_X +l A)\cdot C_j.$$
Now that 
$K_X +l A$
is nef, it follows that
$0=(K_X +l A)\cdot\gamma^+ \geq A\cdot\gamma^+$,
and thus
$\gamma^+ \approx0$.
Hence, we obtain
$$\gamma\approx\sum_{j=1}^r a_j C_j.$$
From (\ref{def}.5),
$$0=H\cdot\gamma +(f^* {\cal A})\cdot\gamma.$$
Since
$H$
and 
$f^* {\cal A}$
are both nef, we have
$$0=H\cdot\gamma=\sum_{j=1}^r a_j H\cdot C_j,$$
and 
$a_j H\cdot C_j =0$
for all $j=1,\ldots, r$.
If there were more than one $j$ such that 
$a_j >0$,
then
$H\cdot C_j =0$
and
$C_j \in R$
for more than one $C_j$,
a contradiction to the choice that 
${C_1},\ldots,{C_r}$
generate distinct extremal rays.
Hence, there is only one non-zero 
$a_{j_0}$ 
and all the others are zero, i.e., 
$\gamma\approx a_{j_0} C_{j_0} \in R$.
 
\end{pf}


\end{document}